\documentclass[12pt]{article}
\usepackage[cm]{fullpage}
\usepackage[small]{titlesec}
\usepackage[cm]{fullpage}

\usepackage{amsmath,amscd, float,rotating}
\usepackage{pb-diagram}

\usepackage{latexsym}
\usepackage{amsfonts}
\usepackage{amsmath}
\usepackage{amssymb}

\numberwithin{equation}{section}

\def\XXint#1#2#3{{\setbox0=\hbox{$#1{#2#3}{\int}$}
\vcenter{\hbox{$#2#3$}}\kern-.5\wd0}}

\newcommand{\dbar}{\overline{\partial}}

\newcommand{\ddt}[1]{\frac{\partial #1}{\partial t}}

\newcommand{\ddbar}{\frac{\sqrt{-1}}{2\pi} \partial\dbar}

\begin{document}
\newcounter{remark}
\newcounter{theor}
\setcounter{remark}{0} \setcounter{theor}{1}
\newtheorem{claim}{Claim}
\newtheorem{theorem}{Theorem}[section]
\newtheorem{proposition}{Proposition}[section]
\newtheorem{lemma}{Lemma}[section]
\newtheorem{definition}{Definition}[section]
\newtheorem{conjecture}{Conjecture}[section]
\newtheorem{corollary}{Corollary}[section]
\newenvironment{proof}[1][Proof]{\begin{trivlist}
\item[\hskip \labelsep {\bfseries #1}]}{\end{trivlist}}
\newenvironment{remark}[1][Remark]{\addtocounter{remark}{1} \begin{trivlist}
\item[\hskip \labelsep {\bfseries #1
\thesection.\theremark}]}{\end{trivlist}}
\newenvironment{example}[1][Example]{\addtocounter{remark}{1} \begin{trivlist}
\item[\hskip \labelsep {\bfseries #1
\thesection.\theremark}]}{\end{trivlist}}
~

\def\b{{\bf \beta}}

\centerline{ \bf \Large The
Ricci flow on the sphere with marked points
\footnote{Work supported in part by
National Science Foundation grants DMS-12-66033, DMS-0847524 and DMS-0905873 and a Collaboration Grants for Mathematicians from Simons Foundation.} }
\bigskip\bigskip

\centerline{D.H. Phong*, Jian Song$^\dagger$,  Jacob Sturm$^\ddagger$ and Xiaowei Wang $^\ddagger$}

\bigskip

\smallskip

\medskip

\begin{abstract}

{\small The Ricci flow on the $2$-sphere with marked points is
shown to converge in all three stable, semi-stable, and unstable cases. In the stable case, the flow was known to converge without
any reparametrization, and a new proof of this fact is given. The semi-stable and unstable cases are new, and it is shown that the flow converges in the Gromov-Hausdorff topology to a limiting metric space which is also a $2$-sphere, but with different
marked points and hence a different complex structure. The limiting metric space carries a unique conical constant curvature metric  in the semi-stable case, and a unique conical shrinking gradient Ricci soliton in the unstable case.}

\end{abstract}

\section{Introduction}

A central theme in geometry is to characterize geometric structures by canonical metrics. The Uniformization
Theorem achieves this for smooth compact Riemann surfaces. In higher dimensions, a well-known conjecture of Yau \cite{Y},
broadly stated, is that the existence of a canonical metric should be equivalent to a suitable notion of stability
in geometric invariant theory. When the structure is not stable, it is expected that a canonical metric should still exist,
albeit with singularities or on an adjacent structure. A model scenario is that of holomorphic vector bundles $E\to M$
over a compact K\"ahler manifold $M$. When $E$ is stable, a Hermitian-Einstein metric will exist, by the celebrated
theorem of Donaldson-Uhlenbeck-Yau \cite{D87,UY}. When $E$ is unstable, the Yang-Mills flow will converge instead
to a Yang-Mills
connection on the double dual of the Harder-Narasimhan-Seshadri filtration of $E$. This last statement was
conjectured by Bando and Siu \cite{BS}. It was proved in ${\rm dim}\,M=2$ by Daskalopoulos and Wentworth \cite{DW},
and very recently in general by Jacob \cite{J}, Sibley \cite{Si}, and Sibley and Wentworth \cite{SiW},
building on the ideas of \cite{D87, UY} and Uhlenbeck \cite{U}.

\medskip

When we pass from holomorphic vector bundles to complex manifolds,
the Yang-Mills flow is replaced by the Ricci flow. The major questions are then to determine the metric,
the complex structure, and the singularities which would emerge from its long-time limit. This question is significantly more difficult
than in the Yang-Mills case, because the equations are more non-linear, and the group of diffeomorphisms
is more subtle than the group of gauge transformations.
The case of positive curvature has proved in particular to be quite challenging,
and there are still few reasonably complete results. Two very recent important advances are the work of Tian-Zhang \cite{TZ} in 3 dimensions,
and the work of Chen-Wang \cite{CWb} in all dimensions, on the convergence of the flow to a soliton with mild singularities. Besides these, there is also
a body of works on Riemann surfaces with marked points which we discuss next in greater detail.

\medskip

First, we need the notion of metrics with conical singularities on a Riemann surface. Let $M$ be a compact Riemann surface, and $p$ a given point on $M$. A metric $g$ on $M$ is said to have a conical singularity at $p$ if it can be expressed as
\begin{eqnarray}
\label{conic}
g= e^{f(z)} |z|^{-2\beta} |dz|^2
\end{eqnarray}
near $p$, with $f(z)$ a bounded function. Here $z$ is a local holomorphic coordinate
centered at $p$, and $\beta\in (0,1)$ is a constant. The constant $\beta$ is sometimes referred to as the weight of $g$ at $p$, and the cone angle of $g$ at $p$ is $2(1-\beta)\pi$. To lighten the notation, we denote by $g$ both the metric and the corresponding K\"ahler form ${\sqrt{-1}\over 2\pi}e^{f(z)}|z|^{-2\beta}|dz|^2$, when which is intended is clear from the context.

\smallskip
More generally, we consider a compact Riemann surface $M$ with given points $p_1,\cdots,p_k$, and weights $\beta_j$ associated to each point $p_j$. We denote by $\beta$ the divisor $\beta=\sum_{j=1}^k\beta_jp_j$, and refer to the data $(M,\beta)$ as a pair. If $g$ is a $C^2$ metric on $M\setminus\beta$, the Ricci curvature can be defined on $M\setminus\beta$ by the usual formula
$$
Ric(g)=-{\sqrt{-1}\over 2\pi}\partial\bar\partial \log\,g
$$
If in addition $g$ admits a conical singularity with weight $\beta_j$ at each point $p_j$, then the Ricci curvature ${\bf Ric}(g)$ can be defined on the whole surface $M$ as a current
$$
{\bf Ric}(g)=Ric(g)
+
\sum_{j=1}^k\beta_j[p_j]
$$
where $[p_j]$ is the Dirac measure at $p_j$. This is because, in two dimensions, the contribution $i\partial\bar\partial f$ of the conformal factor $f$ in (\ref{conic}) cannot include a singular measure  if $f$ is bounded and $C^2$ away from the point $p$.
We restrict ourselves to metrics $g$ with conical singularities whose Ricci current ${\bf Ric} (g)$ is still in $c_1(M)$, i.e., the same Chern class as the Ricci curvature of smooth metrics on $M$. This means that
$$
\chi(M)=\int_M {\bf Ric}(g)=\int_{M\setminus\beta}
Ric(g)+\sum_{j=1}^k\beta_j,
$$
and hence
$$
\int_{M\setminus\beta}
Ric(g)=\chi(M)-\sum_{j=1}^k\beta_j
\equiv\chi(M,\beta).
$$
where the second equality just defines
the Euler characteristic $\chi(M,\beta)$ of the pair $(M,\beta)$.

\smallskip

The metric $g$ with conical singularities is said to have constant Ricci curvature if it satisfies $Ric(g)=\mu\,g$
on $M\setminus\beta$ for some constant $\mu$. In view of the requirement that ${\bf Ric}(g)$ is still in the same Chern class $c_1(M)$, it follows that the constant $\mu$ must satisfy the constraint $\chi(M,\beta)=\mu\,\int_M g$. If we normalize the metric $g$ so that $\int_M g=2$, the equation of constant Ricci curvature becomes
\begin{equation}
\label{ricci}
Ric(g)={1\over 2}\chi(M,\beta) g
\end{equation}
on $M\setminus\beta$.
The metrics with conical singularities and constant Ricci curvature on Riemann surfaces $(M,\beta)$ with weights have been extensively studied by Troyanov \cite{Tr},
Luo-Tian \cite{LT}. When the pair $(M,\beta)$ has Euler characteristic $\chi(M,\b)\leq 0$, it has been shown in \cite{Tr} that it always admits a metric with conical singularities and constant Ricci curvature,
and that such a metric is unique up to scaling. Clearly, $\chi(M,\b)$ can be strictly positive only when
\begin{equation}
\label{pos}
M=S^2,
\qquad
\qquad
\sum_{j=1}^k\beta_j<2.
\end{equation}
Henceforth, we shall make these assumptions. In this case, there are indeed in general obstructions for the existence of metrics with conical singularities and constant Ricci curvature. More specifically, the following is known:

\smallskip

$\bullet$ When $k=1$, equation (\ref{ricci}) does not admit a solution. Instead, one can construct a unique rotationally symmetric compact shrinking soliton $g\in c_1(S^2)$ (the tear drop) \cite{Wu, BM, Ra}.

$\bullet$ When $k=2$, if  $\beta_1 = \beta_2$, there exists a unique rotationally symmetric solution of equation (\ref{ric}) (the football) \cite{Ch, BM, Ra}.

$\bullet$ When $k=2$ and $\beta_1\neq \beta_2$, equation (\ref{ricci}) does not admit a solution. Instead, one can construct a unique rotationally symmetric compact shrinking soliton $g$ \cite{Wu,  BM, Ra}.

$\bullet$ When $k\geq 3$, there does not exist any holomorphic vector field on $S^2$ fixing $p_1, ..., p_k$ since any holomorphic vector field can at most vanish at 2 distinct points. Then the equation (\ref{ricci}) admits a unique solution if and only if \cite{Tr,LT}
 $$
 2\,{\rm max}_j\,\beta_j<\sum_{j=1}^k\beta_j.
 $$

 \medskip
Next, we describe what is known about the Ricci flow on the sphere $S^2$. The case without marked points and conic singularities has been completely settled by the work of Hamilton \cite{H} and Chow \cite{Ch}. The orbifold shrinking gradient solitons on $S^2$ have been classified by Wu \cite{Wu}. The Ricci flow for metrics with conical singularities on Riemann surfaces was first studied
by Yin \cite{Yi1, Yi2}, who provided an important analytic framework as well as a proof of the long-time existence, and convergence of the flow when $\chi(M,\beta)\leq 0$.
Another approach to existence results for the Ricci flow for metrics with conic singularities was given by Mazzeo, Rubinstein and Sesum \cite{MRS}, using an extensive machinery of polyhomogeneous expansions, conormal distributions and b-spaces. They also prove  the convergence of the flow on any stable pair $(S^2, \beta)$. When the pair $(S^2,\beta)$ is not stable, they argue for some notion of ``geometric" convergence (see Theorem 1.3 as well as section \S 5, and especially Proposition 5.3 in \cite{MRS}). The fact is that the case of $(S^2,\beta)$ not stable presents some significant new difficulties which were absent in the stable case and which cannot be bypassed. On one hand, the existence and global structure of the limiting space have to be determined, and it is far from evident a priori that the limiting space is another pair $(S^2,\beta_\infty)$, or even that its singular set is closed. In fact, all conical singularities except the main one may converge to a single limiting conical point, and thus the injectivity radius will converge generically to $0$, preventing any application of Hamilton's compactness theorem. Furthermore, for the possible convergence to a soliton, one must allow for reparametrizations, thus ruling out techniques based solely on multiplier ideal sheaves.

\medskip
The case of pairs $(S^2,\beta)$ is of particular importance as it is the most basic example of spaces that can exhibit all three types of geometric structures, namely stable, semi-stable, and unstable structures. It is also a good model case of the Ricci flow on complex manifolds with singularities. The goal of the present paper is to provide a complete understanding of the long-time behavior of the Ricci flow in this case. Besides its own interest, such an understanding should be valuable in the development of any program to produce canonical metrics on singular K\"ahler manifolds by the Ricci flow.

\medskip
We state now precisely our results. As in \cite{Yi1, Yi2}, we fix a metric $g_\beta$ with conic singularities on the pair $(S^2,\beta)$, which is smooth away from the points $\{p_1,\cdots,p_k\}$, given by $g_\beta= |z|^{-2\beta_i}|dz|^2$ in a holomorphic coordinate $z$ near each point $p_i$,
and normalized to that $\int_{S^2}g_\beta=2$. Note that this normalization coincides with $g\in c_1(S^2)$.
The Ricci flow is the following flow of metrics $g(t)$,
\begin{equation}
\label{ric}
{\partial g(t)
\over
\partial t}=-Ric(g)+{1\over 2}\chi(S^2,\beta)g(t),
\qquad
g(0)=g_0
\end{equation}
on $S^2\setminus\beta$, where $g_0$ is a given initial metric, also normalized to have area $2$, $\int_{S^2}dg_0=2$. Since the Ricci flow preserves the conformal class of the metric,
we can write $g(t)=e^{u(t)}g_\beta$, and the Ricci  flow is equivalent to
the following flow for the conformal factor $u(t)$,
\begin{equation}
\label{ricu}
{\partial u(t)\over\partial t}
=
e^{-u(t)}\Delta_{g_\beta}u(t)
+
{1\over 2}\chi(S^2,\beta)
-e^{-u(t)}R_{g_\beta},
\qquad
u(0)=u_0
\end{equation}
where $R_g$ and $\Delta_g$ refer to the curvature of the Laplacian of a given metric $g$, and we have expressed the initial metric $g_0$ as
\begin{equation}
\label{u0}
g_0= e^{u_0}g_\beta.
\end{equation}

Because of the conical singularities, it is important to specify the regularity of the metrics $g(t)$, $t\geq 0$. In \cite{Yi1, Yi2}, the key notion of weighted Schauder spaces $C^{\ell,\alpha}(S^2,\beta)$ for a surface with conical singularities was introduced, for each $\ell\in {\bf N}$, $\alpha\in (0,1)$. Their precise definition will be recalled in Section \S 2 below.
The following
conditions were imposed on
the initial metric $g_0=e^{u_0}g_\beta$,
\begin{eqnarray}
\label{Yinconditions}
&&
u_0\in C^{2,\alpha}(S^2,\beta), \quad \int_{S^2}|\nabla u_0|^2 dg_0<\infty,
\nonumber\\
&&
R_{g_0}\in C^{2,\alpha}(S^2,\beta), \quad \int_{S^2}|\nabla R_{g_0}|^2 dg_0<\infty,
\quad
\Delta_{g_0}R_{g_0}\in L^\infty(S^2)
\end{eqnarray}
Here $R_{g_0}$ is the curvature of $g_0$, and $\Delta_{g_0}$ is the Laplacian with respect to $g_0$. It was then shown by Yin \cite{Yi2} that, under the conditions (\ref{Yinconditions}) on the initial metric $g_0$,
the Ricci flow will then exist for all time, and the conditions (\ref{Yinconditions}) are preserved  for all $t\geq 0$. Furthermore, the area normalization and the total curvature are also preserved,
\begin{equation}
\label{preserve}
\int_{S^2} dg(t)=2,
\qquad
\int_{S^2\setminus\beta}R(t) dg(t)
=\chi(S^2,\beta).
\end{equation}

\medskip
For our purposes, we shall consider a slightly more regular class of metrics which can be defined as follows. Let ${\cal B}(S^2,\beta)$ be the following space
\begin{equation}
\label{calB}
{\cal B}(S^2,\beta)
=
C^{2,\alpha}(S^2,\beta) \cap W^{1,2}(S^2).
\end{equation}
Here $W^{1,2}(S^2)$ is the usual Sobolev space with respect to the metric $g_\beta$. But since the Dirichlet energy $\int_{S^2}|\nabla u|_\beta^2dg_\beta$ is conformally invariant in 2 dimensions, and since the space $C^{2,\alpha}(S^2,\beta)$ is already included in $L^\infty(S^2)$, the space ${\cal B}(S^2,\beta)$ is actually independent of the choice of metric within the conformal class.
We can now introduce the class of metrics that we shall consider:

\begin{definition}
\label{regularA}

A metric $g=e^{u}g_\beta$ is said to be a {\it regular metric with conical singularities} if
\begin{equation}
\label{regularB}
u\in {\cal B}(S^2,\beta),
\quad
R_g\in {\cal B}(S^2,\beta),
\quad
\Delta_g R_g\in {\cal B}(S^2,\beta).
\end{equation}

\end{definition}

\medskip
As we shall see, the stronger requirement that $\Delta_g R_g\in {\cal B}(S^2,\beta)$ in Definition \ref{regularA} compared to Yin's requirement that $\Delta_g R_g\in L^\infty(S^2)$ as in (\ref{Yinconditions}) is what allows an extension of Perelman's methods \cite{P} to the Ricci flow with regular initial metric.

\medskip

It is also convenient to introduce the following terminology:

\medskip

\begin{definition} Let $(S^2, \beta)$ be a sphere with marked points, with $k\geq 3$.
We shall say that

$\bullet$ $(S^2, \beta)$ is stable if $\sum_{i=1}^k \beta_i \geq 2$ or $2\,{\rm max}_j\beta_j < \sum_{i=1}^k \beta_i$.

$\bullet$ $(S^2, \beta)$ is semi-stable if $\sum_{i=1}^k \beta_i <2$ and $2\,{\rm max}_j\beta_j = \sum_{i=1}^k \beta_i$.

$\bullet$ $(S^2, \beta)$ is unstable if $\sum_{i=1}^k \beta_i <2$ and $2\,{\rm max}_j\beta_j > \sum_{i=1}^k \beta_i$.

\end{definition}

 Without loss of generality, we assume that
 $$ \beta_1 \leq \beta_2 \leq ... \leq \beta_k= \beta_{max}. $$
 We have then the following

\medskip
\begin{theorem}\label{main}  Let $(S^2,\beta)$
be a sphere with $k$ marked points,
$\beta=\sum_{j=1}^k\beta_jp_j$,
$\sum_{j=1}^k\beta_j<2$ and $k\geq 3$.
Consider the Ricci flow $g(t)$ with an initial metric $g_0$ which is regular in the sense of Definition \ref{regularA}, and normalized to $\int_{S^2}dg_0=2$. Then the metrics $g(t)$ are also regular in the sense of Definition \ref{regularA} for each $t\geq 0$. Furthermore, the behavior of the flow as $t\to+\infty$ can be described as follows:

\begin{enumerate}

\item If $(S^2,\beta)$ is stable, then the flow converges in the Gromov-Hausdorff topology and
in $C^{\infty}(S^2\setminus\beta)$ to the unique conical constant curvature  metric $g_\infty\in c_1(S^2)$  on $(S^2,\beta)$.

\item If $(S^2,\beta)$ is semi-stable, then the flow converges in the Gromov-Hausdorff topology to the unique conical constant  curvature metric $g_\infty$
on a pair $(S^2,\beta_\infty)$, where the divisor $\beta_\infty$ is given by
$$
\beta_\infty= \beta_{max} [p_\infty] + \beta_{max} [q_\infty],
$$
with $\beta_{max}\equiv {\rm max}_{1\leq j\leq k}\beta_j$.
The   convergence is smooth on  $S^2\setminus \{p_\infty, q_\infty\}$. In particular,  $p_k$ converges in Gromov-Hausdorff distance to one of the two points $p_\infty$ and $q_\infty$, while  $p_1, ..., p_{k-1}$  converge to the other.

\item If $(S^2,\beta)$ is unstable, then the flow converges in the Gromov-Hausdorff topology to the unique conical  shrinking gradient Ricci soliton $g_\infty$ on a pair $(S^2, \beta_\infty)$ with
$$\beta_\infty = \beta_{p_\infty} [p_\infty] + \beta_{q_\infty} [q_\infty]~,  0\leq \beta_{q_\infty} < \beta_{p_\infty}, ~ \beta_{p_\infty} + \beta_{q_\infty} = \sum_{i=1}^k \beta_i. $$
Furthermore, the limiting space $(S^2, \beta_\infty)$ and limiting metric $g_\infty$ do not depend on the initial metric $g_0$.

\end{enumerate}

\end{theorem}

\bigskip

We remark that the stable case (1) in Theorem \ref{main} is proved in \cite{MRS}.

\bigskip

Under an additional assumption, we can determine the limiting points $p_\infty$ and $q_\infty$ in the unstable case as well:

\begin{theorem} \label{main2}
Assume that $(S^2, \beta)$ is unstable, so that we are in the case 3)
of Theorem \ref{main}. There exists an explicit constant $  \mathcal{W}_\beta$ such that if Perelman's conical W-functional is strictly bounded below by $\mathcal{W}_\beta$ at $t=0$ for some initial regular metric $g_0$,
\begin{equation}\label{wlo}
\mu(g_0) > \mathcal{W}_\beta,
\end{equation}
then $p_k$ converges to $p_\infty$ and $p_1,\cdots,p_{k-1}$ converge to $q_\infty$. In particular,
$$
\beta_\infty = \beta_{max} [p_\infty] + (\sum_{i<k}  \beta_i ) [q_\infty].
$$
\end{theorem}

\medskip
The constant $\mathcal{W}_\beta$ can be calculated explicitly (see Proposition \ref{soon} below). However, we believe that one should always be able to construct an initial metric $g_0$ satisfying the condition in Theorem  \ref{main2} by perturbation and gluing, so that the conclusion of Theorem \ref{main2} should always hold.

\bigskip

The paper is organized as follows. In Section \S 2, we recall Yin's results \cite{Yi1,Yi2} on the linear heat equation, the long-time existence of the Ricci flow, as well as his notions of weighted Schauder spaces for conic pairs $(S^2,\beta)$. We deduce that the Ricci flow preserves the class of ``regular metrics" introduced in Definition \ref{regularA}. In Section \S 3, we provide an equivalent formulation of the Ricci flow as a flow of K\"ahler potentials. Section \S 4 is a first major step: it shows how, by restricting the Ricci flow to the class of regular metrics introduced in Definition \ref{regularA}, the Perelman monotonicity can be extended to the case of Riemann surfaces with conic singularities. As a consequence, we obtain the essential fact that the scalar curvature and the diameter are uniformly bounded along the flow. Since in two dimensions, the scalar curvature determines the full Riemannian curvature, we are in a situation similar to that considered in \cite{PS06}, with the key additional complication that the manifolds are not compact, and the injectivity radius not bounded from below. Section \S 5 is devoted to the analysis of the long-time behavior in the stable case. In this case, the functional $F_\beta$ is proper, and we adapt the arguments in the smooth case to show the convergence of the flow. Section \S 6 is devoted to the semi-stable and the unstable cases. We combine Perelman monotonicity, Cheeger-Colding theory, the partial $C^0$ estimate, and Hamilton's entropy to establish the sequential convergence of the flow, in the Gromov-Hausdorff sense, to a sphere with marked points, equipped with either a metric of constant curvature or a shrinking gradient soliton. In the semi-stable case, we derive a weak lower bound for the functional $F_\beta$ which allows us nevertheless to show the existence of a sequence of times along which the Ricci flow converges to a metric of constant curvature. In the unstable case, we show that the limiting metric cannot have constant curvature, and hence it must be a soliton. Finally, Section \S 7 is devoted to the proof of Theorem \ref{main2}. A key tool is Perelman's functional $W(g,f,\tau)$, with $\tau$ chosen to be the singular time of the flow.


\section{The Ricci flow with initial regular metric}

The purpose of this section is to show that the Ricci flow preserves the regularity of the initial metric.

\subsection{Results of Yin on the heat equation and the Ricci flow}

We begin by recalling some results of Yin \cite{Yi1, Yi2} on the heat equation and the Ricci flow on Riemann surfaces with conical singularities which play an essential role in the sequel.

\smallskip
First, we recall the definition of the weighted Schauder norms introduced in \cite{Yi1, Yi2}. Fix a conical singularity $p$ with weight $\beta$ in a pair $(S^2,\beta)$, as described in section \S 1. Let $z=re^{i\theta}$ be a complex coordinate centered at $p$, defined in a disk $U$. Then the $C^{\ell,\alpha}(U,p)$ norm of a function $f$ on $U$ is defined by
\begin{equation}
\label{Holder}
\|f\|_{C^{\ell,\alpha}(U,p)}
=
{\rm sup}_{m\in {\bf N}}\|F_m(s,\theta)\|_{C^{\ell,\alpha}((2^{-1}\leq s\leq 2)\times S^1)}
\end{equation}
where $F_m(s,\theta)\equiv f((1-\beta)^{1\over 1-\beta}(2^{-m}s)^{1\over 1-\beta} e^{i\theta})$, and the $C^{\ell,\alpha}$ norm on the right hand side is the usual Schauder norm on the annulus $(2^{-1}\leq s\leq 2)\times S^1$. The weighted Schauder norm of a function on a pair $(S^2,\beta)$ is obtained by covering $S^2$ by a finite number of neighborhoods $V_\kappa$, $V_\kappa\subset\subset S^2\setminus\beta$, together with disks $D_j$ centered at each conical singularity as described above. We can then set
\begin{equation}
\|f\|_{C^{\ell,\alpha}(S^2,\beta)}
=
\sum_\kappa\|f\|_{C^{\ell,\alpha}(V_\kappa)}
+
\sum_{j=1}^k
\|f\|_{C^{\ell,\alpha}(U_j,p_j)}.
\end{equation}
Similarly, we can define weighted parabolic Schauder norms $C^{k,\alpha}((S^2,\beta)\times [0,T])$ for functions $f(z,t)$ defined on $(S^2,\beta)\times [0,T]$, with the weighted norm $C^{\ell,\alpha}((U,p)\times [0,T])$ defined now by
\begin{equation}
\|f\|_{C^{\ell,\alpha}((U,p)\times [0,T])}
=
{\rm sup}_{m\in {\bf N}}\|F_m(s,\theta,\tilde t)\|_{C^{\ell,\alpha}
((2^{-1}\leq s<2)\times S^1\times [0,2^{2m}T])}
\end{equation}
and $F_m(s,\theta,\tilde t)
\equiv  f((1-\beta)^{1\over 1-\beta}(2^{-m}s)^{1\over 1-\beta} e^{i\theta}, 2^{-2m}\tilde t)$. The global weighted Schauder norms on $(S^2,\beta)\times [0,T])$ can then be defined by patching up as above.

\medskip
The theorems of Yin that we need are the following: the first is the following theorem on the existence and estimates for the Ricci flow for conical metrics:

\begin{theorem}
\label{Yin1}
(\cite{Yi2}, Theorem 4.1 and end of Section \S 5.1) Let $g_0=e^{u_0}g_\beta$ be a metric satisfying the conditions (\ref{Yinconditions}).
Then the Ricci flow (\ref{ric}) exists for all time $t\in [0,\infty)$, and the conditions
$\int_{S^2}dg=2$, $\int_{S^2\setminus\beta}R=\chi(S^2,\beta)$ are both preserved. For each $t\in [0,\infty)$,
the metric $g(t)=e^{u(z,t)}g_\beta$ again satisfies the conditions (\ref{Yinconditions}).

Furthermore, the conformal factor $u$ and the curvature $R$ of $g(t)$ satisfy natural parabolic estimates, in the sense that for each fixed $T<\infty$, the conformal factor $u$ is in $C^{2,\alpha}((S^2,\beta)\times [0,T])$, and the curvature $R$ is in $C^{2,\alpha}((S^2,\beta)\times [0,T])$.

\end{theorem}

\noindent
From the general theory of parabolic equations, it is also not hard to see that for $t>0$, the metric $g(t)$ is smooth away from the conical singularities. As mentioned earlier, the area and the Chern class are preserved by the flow (\ref{preserve}).
We observe that the higher regularity $R\in C^{2,\alpha}((S^2,\beta)\times [0,T])$ stated above follows from the second result of Yin that we need,
namely the following existence, regularity, and uniqueness theorem on the heat equation on weighted H\"older spaces, where we have introduced the following parabolic analogue of the space ${\cal B}(S^2,\beta)$,
$$
{\cal B}((S^2,\beta)\times[0,T])\equiv
 C^{0, \alpha}((S, \beta)\times [0, T])\cap C^0([0,t],W^{1,2}).
$$

\begin{theorem}
\label{Yin2}
(\cite{Yi2}, Lemmas 2.2 and 2.4) Consider the following linear heat equation
\begin{eqnarray}
\label{linear}
{\partial u\over\partial t}
=
a(z,t)\Delta_{g_\beta}u +b(z,t)u +f(z,t),
\qquad
u(z,0)=u_0(z).
\end{eqnarray}
for a function $u$ on $(S^2,\beta)\times [0,T]$, where $a(z,t)>0$.

{\rm (1)} If $a, b, f\in C^{0, \alpha}((S^2, \beta)\times [0, T])$,  $u_0\in C^{2, \alpha}(S^2, \beta)$, then the equation admits a unique solution $u(z, t)\in C^{2, \alpha}((S^2, \beta)\times [0, T])$.

{\rm (2)}  If $a,b,f\in
 {\cal B}((S^2,\beta)\times[0,T])$ and $u(\cdot,0)\in {\cal B}(S^2,\beta)$, then
 $u\in {\cal B}((S^2,\beta)\times[0,T])$.

{\rm (3)}  If $\partial_ta,\partial_tb,\partial_tf\in
 {\cal B}((S^2,\beta)\times[0,T])$ and $\partial_tu(\cdot,0)\in {\cal B}(S^2,\beta)$, then
 $\partial_tu\in {\cal B}((S^2,\beta)\times[0,T])$.

\end{theorem}

We observe that the statement (3) in the above theorem follows from the fact that $u={\rm lim}\,u_k$, where $u_k$ satisfies the heat equation, with Neumann condition, on the complement
of the union of disks centered at each marked point. The functions $\partial_tu_k$ satisfy the same Neumann condition, and arguing as in \cite{Yi1}, Theorem 3.1, shows that a subsequence converges to a function $v\in C^{2,\alpha}((S^2,\beta)\times[0,T])$. Since $\partial_tu_k$ converges weakly to $\partial_tu$, we must have $v=\partial_tu$. The fact that $\partial_tu\in
{\cal B}((S^2,\beta)\times[0,T])$ follows then from (2) of Theorem \ref{Yin2}.

\subsection{Preservation of regularity}

We establish now the first statement in Theorem \ref{main}, namely that regularity in the sense of Definition \ref{regularA} is preserved along the Ricci flow.
Assume then that the initial metric $g_0=e^{u_0}g_\beta$ is regular in the sense of Definition \ref{regularA}. In particular, it satisfies the conditions of Yin's Theorem \ref{Yin1}, and it suffices to show that if $R(t)$ denotes the curvature of $g(t)$
and $\Delta_t$ denotes the Laplacian with respect to $g(t)$, then $\Delta_t R(t)$ is in ${\cal B}((S^2,\beta)\times [0,T])$.

Now it is well-known that, under the Ricci flow, the curvature flows by
\begin{equation}
\label{Rflow}
\partial_t R
=
\Delta_tR+(R-{1\over 2}\chi(S^2,\beta))R.
\end{equation}
Differentiating the equation gives the flow for $\partial_tR$,
\begin{equation}
\label{Rtflow}
\partial_t(\partial_tR)
=
\Delta_t(\partial_tR)
+R\partial_t R+R(R-{1\over 2}\chi(S^2,\beta))^2
\end{equation}
where we also made use of the defining equation for the Ricci flow,
$\partial_tu=-R+{1\over 2}\chi(S^2,\beta)$.
Note that, in general, the partial derivative $\partial_tu$ of the solution of a heat equation such as (\ref{linear}) satisfies another heat equation, but whose coefficients are less regular. The basic observation is that, in the case of $R(t)$ in dimension 2,
its derivative $\partial_t R(t)$ satisfies a heat equation
whose coefficients are as regular as those of the heat equation for $R(t)$. Now at time $t=0$,
\begin{equation}
\partial_tR=\Delta_0R_0+(R_0-{1\over 2}\chi(S^2,\beta))R_0
\in {\cal B}(S^2,\beta),
\end{equation}
since both $R_0$ and $\Delta_0R_0$ are in ${\cal B}(S^2,\beta)$ by hypothesis.
By Yin's Theorem \ref{Yin2}, the equation for $\partial_tR$ admits a unique solution in ${\cal B}((S^2,\beta)\times [0,T])$. To show that this solution coincides with $\partial_tR$, we can argue as in the proof of the statement (3) of Theorem \ref{Yin2}),
by considering the limits of $\partial_tR_k$, where $R_k$ are the solutions of the corresponding Neumann problems. Thus $\partial_tR(t)\in {\cal B}((S^2,\beta)\times[0,T])$. Using the heat equation (\ref{Rflow}) for $R(t)$, we deduce that $\Delta_tR(t)\in {\cal B}((S^2,\beta)\times[0,T])$, as was to be proved.

\subsection{Continuity of the conformal factor $u(t)$ at the marked points
$p_1,\cdots,p_k$}

In general, a function with finite $C^{\ell,\alpha}(S^2,\beta)$ has bounded $L^\infty(S^2)$ norm, but may not be necessarily continuous at the conical points $p_1,\cdots,p_k$. We observe however that, as a consequence of the boundedness of the curvature on each finite-time interval, the conformal factors $u(t)$ in the Ricci flow are always continuous:

\begin{lemma} Suppose $g(t) = e^{u(t)}g_\beta$ is a solution in $t\in[0,T)$ of the Ricci flow with a regular initial metric $g_0$. Then $u(t)\in C^0(S^2\times [0, T))$.

\end{lemma}

\noindent
{\it Proof}. We already know that $u(t)$ and $R(t)$ are in $L^\infty(S^2\times [0, T'])$ for any $T'<T$. Since
\begin{equation}
\Delta_t u(t) = e^{-u(t)} R_{g_\beta}- R(t).
\end{equation}
it follows that $\Delta_t u(t)$ is uniformly bounded. By Theorem 3.1 and Lemma 3.2 in \cite{Yi2},
$| \nabla u(t) |_{g(t)}$
is bounded. Since $g(t)$ and $g_\beta$ are equivalent in $L^\infty$,
$| \nabla u(t) |_{g_\beta}$ is also bounded.
Now for a fixed conical point $p$ and any two points $q_1$, $q_2 \in B_{g_\beta} (p, 1)$,
$$|u(t, q_1)- u(t, q_2)|  \leq \sup_{x\in S^2} |\nabla  u(t, x)|_{g_\beta} d_{g_\beta} (q_1, q_2) \leq K d_{g(t)} (q_1, q_2), $$
where $d_{g_\beta}(q_1, q_2)$ and $d_{g(t)}(q_1, q_2)$ are the geodesic distance between $q_1$ and $q_2$ with respect to $g_\beta$ and $g(t)$ respectively. Thus $u(t,q)$ is Lipschitz in $q$, uniformly in $t$. It is also Lipschitz in $t$, since $\partial_t{u} = e^{-u} \Delta_{g_\beta} u + {1\over 2}
\chi(S^2,\beta)- e^{-u} R_{g_\beta}$ and hence $\partial_tu$ is uniformly  bounded.


\section{The point of view of complex geometry}

So far we have expressed the Ricci flow and the metrics $g(t)$ in terms of their conformal factors $g(t)=e^{u(t)}g_\beta$ and in terms of real differential geometry. But since the long-time behavior of the flow will ultimately depend on the complex structure of the pair $(S^2,\beta)$, we have to formulate the flow and the metrics in terms of K\"ahler potentials and in terms of complex differential geometry. This formulation is analogous to the smooth case, with suitable modifications required by the presence of conical singularities.

\subsection{K\"ahler potentials and Ricci potentials}

Fix a metric $g_\beta$ satisfying the two conditions
\begin{equation}
\label{normalization}
\int_{S^2} dg_\beta=2,\qquad\int_{S^2\setminus\beta} R_\beta=\chi(S^2,\beta).
\end{equation}
If $g$ is a regular metric satisfying the first condition, then by
Yin's Theorem 3.1 and Lemma 3.2 \cite{Yi2}, we can write
\begin{equation}
\label{kahlerpotential}
g=g_\beta+{\sqrt{-1}\over 2\pi}\partial\bar\partial \varphi,
\end{equation}
for some function $\varphi$ which is bounded, has finite Dirichlet energy,
and satisfies the condition that the right hand side be a strictly positive $(1,1)$-form on $S^2\setminus\beta$. Furthermore, $|\nabla \varphi|_\beta^2$ has finite Dirichlet energy. For our purposes, we shall refer to the class of functions satisfying all these properties as $PSH(g_\beta)$. The function $\varphi$ satisfying (\ref{kahlerpotential}) is unique up to a constant.

We define the Ricci potential $h_\beta$ of the metric $g_\beta$ by
\begin{equation}
\label{riccipotential}
Ric(g_\beta)={1\over 2}\chi(S^2,\beta) g_\beta
+
{\sqrt{-1}\over 2\pi}\partial\bar\partial h_\beta,
\qquad
\int_{S^2}e^{h_\beta}dg_\beta=2.
\end{equation}
Again $h$ is bounded, has finite Dirichlet energy, and $|\nabla h|_\beta^2$ has finite Dirichlet energy.

It is now easily verified that, just as in the case of the Ricci flow on smooth K\"ahler manifolds, the Ricci flow is given by $g(t)=g_\beta+{\sqrt{-1}\over 2\pi}\partial\bar\partial\varphi(t)$, with $\varphi(t)$ satisfying the equation
\begin{equation}
\label{rickahler}
\ddt\varphi(t)
=
\log {g_\beta+{\sqrt{-1}\over 2\pi}\partial\bar\partial\varphi(t)\over g_\beta}
+{1\over 2}\chi(S^2,\beta)\varphi -h_\beta.
\end{equation}
The choice of constants in $\varphi(t)$ can be fixed as in \cite{PSS}.

\medskip
It will also be useful to express the Ricci flow with the Fubini-Study metric $g_{FS}$ on $S^2$ instead of $g_\beta$ as reference metric. Note that $g_{FS}=Ric(g_{FS})$, so that the Fubini-Study metric $g_{FS}$ satisfies the first normalization condition in (\ref{normalization}), but not both. Thus we introduce
\begin{equation}
g_{FS,\beta}=g_{FS}^{{1\over 2}\chi(S^2,\beta)}\prod_{j=1}^k|\sigma_j(z)|^{-\beta_j},
\end{equation}
where $\sigma_j(z)$ are holomorphic vector fields with $[2p_j]$ as divisor, and normalized so that $\int_{S^2}dg_{FS,\beta}=2$. The metric $g_{FS,\beta}$ satisfies now both normalization conditions in (\ref{normalization}), and the formulation (\ref{rickahler}) of the Ricci flow applies with $g_\beta$ replaced by $g_{FS,\beta}$. Finally, to rewrite the flow with $g_{FS}$ as reference metric, we set
\begin{equation}
g(t)=g_{FS}+{\sqrt{-1}\over 2\pi}\partial\bar\partial\psi(t)
\end{equation}
Rewriting the Ricci flow with $g_{FS,\beta}$ as reference metric and $\varphi(t)$ as potentials in terms of $g_{FS}$ as reference metric and $\psi(t)$ as potentials, we find
\begin{equation}
\label{ma}
\ddt\psi(t)
=
\log({g_{FS}+{\sqrt{-1}\over 2\pi}\partial\bar\partial\psi\over g_{FS}})
+
{1\over 2}\chi(S^2,\beta)\psi
+
{1\over 2}\sum_{j=1}^k\beta_j\log \,{|\sigma_j|^2
\over g_{FS}}
\end{equation}

\subsection{The $\alpha$-invariant for $(S^2,\beta)$}

The $\alpha$-invariant \cite{T1} can be readily extended to the case of the sphere with marked points. We define
\begin{eqnarray}
\alpha (S^2, \beta) = {\rm sup}\,\alpha
\end{eqnarray}
where $\alpha$ satisfies the condition
\begin{eqnarray}
{\rm sup}_{\varphi \in PSH(g_{FS})}
\int_{S^2} e^{ -{1\over 2}\alpha \chi(S^2,\beta) (\varphi -\sup_{S^2} \varphi) } \prod_{i=1}^k |\sigma_i|^{-\beta_i} g_{FS}^{{1\over 2}\chi(S^2,\beta)} <\infty.
\end{eqnarray}

The following lemma is due to Berman \cite{Be}. We reproduce the short proof below, for the convenience of the reader:

\begin{lemma}\label{alcal}
Assume that $\sum_{i=1}^k \beta_i <2$. Then
\begin{equation}
\alpha(S^2, \beta)= \frac{1-\beta_k}{\chi(S^2,\beta)}.
\end{equation}
\end{lemma}

\medskip
\noindent
{\it Proof.} By a theorem of Demailly \cite{ChS} in the smooth case,
and extended to the conical case by Berman \cite{Be},
the $\alpha$-invariant $\alpha(S^2, \beta)$ is equal to the log canonical threshold
$Lct(S^2,\beta)$, which is defined as follows. Let $h$ be any hermitian metric on $K_{S^2}^{-1}$. Then
\begin{eqnarray}
Lct(S^2, \beta) ={\rm sup}\,\alpha
\end{eqnarray}
where $\alpha$ satisfies the condition
\begin{eqnarray}
\int_{S^2} |\sigma|^{-{\alpha\over m}\chi(S^2,\beta)}  \prod_{i=1}^k |\sigma_i|^{-\mathbb{\beta}_i}
h^{{1\over 2}\chi(S^2,\beta)(1-\alpha)} < \infty.
\end{eqnarray}
for any $m\in \mathbb{Z}^+$, $\sigma \in H^0(S^2, K_{S^2}^{-m}) $.

Since $\beta_k$ is the largest among $\beta_i$, we can calculate the integral near $p_k$. Without loss of generality we may assume $p_k=0$. Then for any $0\leq l \leq 2m$ there is a $\sigma\in H^0(S^2, K_{S^2}^{-m})$ which admits an expansion $\sigma= z^{\ell} f(z) $ near $0$ for some holomorphic function $f(z)$ satisfying $f(0)=1$. Then we have
\begin{eqnarray*}
\int_{S^2} |\sigma|^{ -{ \alpha \over m}\chi^(S^2,\beta)}  \prod_{i=1}^k |\sigma_i|^{-\mathbb{\beta}_i}
&=& \sqrt{-1} \int_{ |z| \leq 1} |z|^{-2{\alpha\ell\over m}
\chi(S^2,\beta)}  |z|^{-2\beta_k} dz\wedge d\bar z + O(1)\\
&\leq& \sqrt{-1} \int_{|z| \leq 1 } |z|^{-2\alpha\chi(S^2,\beta) - 2 \beta_k} dz\wedge d\bar z + O(1)
\end{eqnarray*}
which is finite for any $\alpha< \frac{1-\beta_k}{\chi(S^2,\beta)}$. The equality follows easily by applying a test holomorphic section $\sigma$ in the anti-pluricanonical system.
The proof is complete.

\subsection{The $F$ functional for pairs $(S^2,\beta)$}

It is well-known in the smooth case that the equation of constant scalar curvature admits a variational formulation. As for the $\alpha$-invariant, this can be readily extended to the case of Riemann surfaces with marked points. We can write the corresponding functional in two different ways, depending on whether we use the Fubini-Study metric $g_{FS}$ or the metric $g_\beta$ with conical singularities as reference metric.

With $g_{FS}$ as reference metric, we set
\begin{eqnarray}
\label{FFS}
F_\beta(\varphi) = \frac{\sqrt{-1} }{8\pi} \int_{S^2} \partial \varphi \wedge \bar\partial \varphi   - \frac{1}{2}  \int_{S^2} \varphi ~ dg_{FS}
- {2\over\chi(S^2,\beta)}  \log \left(  \int_{S^2}  e^{ - {1\over 2}
\chi(S^2,\beta) \varphi }  \prod_{i=1}^k |\sigma_i|^{- \beta_i}  g_{FS}^{
{1\over 2}\chi(S^2,\beta)} \right) \nonumber
\end{eqnarray}
while with $g_\beta$ as reference metric, we set
\begin{eqnarray}
\label{Fbeta}
F_\beta(\varphi)
=
{\sqrt{-1}\over 8\pi}\int \partial\varphi\wedge\bar\partial\varphi
-
{1\over 2}\int\varphi dg_\beta
-
{2\over\chi(S^2,\beta)}\log (\int_{S^2}e^{-{1\over 2}\chi(S^2,\beta)\varphi+h_\beta}
dg_\beta).
\end{eqnarray}
In view of the fact that the potentials are always bounded with bounded Dirichlet energy, integration by parts is justified and the two formulations of the $F_\beta$ functional can be verified to agree.
The Euler-Lagrange equation for $F_\beta$ is exactly the equation for the stationary points of the flow (\ref{ma}). It is a special case of the functional $F_\beta$ defined in \cite{SW} for paired Fano manifolds and it satisfies the co-cycle condition.

\medskip
We shall need the following simple property of $F_\beta$, which is a straightforward adaptation of the similar property established in the smooth case in \cite{S}:

\begin{lemma}
\label{proper}
Let $(S^2,\beta)$ be a pair with $\sum_{j=1}^k\beta_j<2$.
 If $\alpha(S^2,  \beta) > 1/2$, then there exists $\epsilon>0$ and $C_\epsilon>0$ such that for all $\varphi \in PSH(S^2, g_{FS})$,
\begin{equation}
 F_\beta(\varphi) \geq \epsilon  \frac{ \sqrt{-1} }{2\pi}   \int_{S^2} \partial \varphi \wedge \bar\partial \varphi  - C_\epsilon
 \end{equation}
In particular, the equation (\ref{ma}) is solvable.

\end{lemma}

We remark that when $(S^2, \beta)$ is not stable, the functional $F_\beta$ is not bounded below and so $(S^2, \beta)$ does not admit a constant curvature metric. Combined with Lemma \ref{alcal}, this can provide a complex geometric proof for the criterion of Troyanov and Luo-Tian as suggested in \cite{Be}.


\section{Perelman monotonicity for regular metrics with conic singularities}

A key result of Perelman \cite{SeT} for the Ricci flow on compact K\"ahler manifolds is the uniform boundedness of the scalar curvature, the uniform boundedness of the diameter, and the uniform lower bound for the volume of balls of fixed center and radii as time evolves. These results follow from the monotonicity of Perelman's famous functional $\mu(g,\tau)$. The purpose of the present section is to extend these results to the Ricci flow on the sphere $(S^2,\beta)$ with marked points, assuming that the initial metric $g_0$ is a regular metric with conical singularities in the sense of Definition \ref{regularA}.

\bigskip
Let $(S^2,\beta)$ be a sphere with marked points, and let $g$ be a regular metric with conical singularities
in the sense of Definition \ref{regularA}. The $W$-functional can be defined for the metric $g$ on $(S, \beta)$ by the same expression as in the smooth case,
\begin{equation}
\label{W}
W(g, f, \tau) = \int_{S^2\setminus \beta } (\tau(R + |\nabla f|^2) + f -2){ e^{-f} \over 4\pi\tau} dg,\qquad\tau>0,
\end{equation}
but the range of the functions $f$ has to be specified with care.
In the smooth case, the function $f$ was simply required to be a smooth function. But in the present case of conical singularities, several requirements have to be simultaneously met: the functional $W(g,f,\tau)$ has to be uniformly bounded from below over the range of $f$; the range of $f$ has to be preserved by the backwards heat flow coupled with the Ricci flow; and the integration by parts required for the monotonicity of Perelman's functional must be justified. We shall show that the following requirements will do the job. Let $F={e^{-f}\over 4\pi\tau}$. Then we set
\begin{equation}
\label{mu}
\mu(g,\tau)={\rm inf}_F W(g,f,\tau)
\end{equation}
where the infimum is taken over all functions $F$ satisfying
\begin{equation}
\label{Fspace}
F\in {\cal B}(S^2,\beta), \quad F\geq \delta_F>0,
\quad
\Delta F\in {\cal B}(S^2,\beta)
\end{equation}
together with the normalization condition
\begin{equation}
\label{Fnormalization}
\int_{S^2}F dg=1.
\end{equation}
Here $\delta_F$ is a strictly positive constant which depends on $F$.

\subsection{The Sobolev inequality on manifolds with conic singularities}

It is well-known since the works of Perelman \cite{P} that the finiteness of the lower bound for $W(g,f,\tau)$ is a consequence of the logarithmic Sobolev inequality, which is itself a straightforward consequence of the Sobolev inequality, by applying Jensen's inequality
\cite{Rot}. Thus the issue is to establish a Sobolev inequality in the case of conical singularities. This can actually be done quite generally in all dimensions, using approximations of K\"ahler metrics with conical singularities by smooth K\"ahler metrics with lower bounds on their Ricci curvatures.

\smallskip

Thus let $(X, L)$ be a polarized K\"ahler manifold of complex dimension $n$. Suppose $D$ is a smooth divisor with $[D] =[L]$. Let $s_D$ be the defining holomorphic section of $D$ and $h$ a smooth positively curved hermitian metric on $L$. A conical K\"ahler metric $\omega$ on $(X,\beta D)$
with conic angle $2\pi(1-\beta)\in (0, 2\pi)$ is a K\"ahler metric on $X\setminus D$
with bounded local potentials, and which is locally equivalent near $D$,
to the standard conical metric $\omega_\beta = \ddbar |z_n|^{2\beta} + \ddbar |z'|^2$, where $\{z_n=0\}$ is a local equation for $D$ and $z'=(z_1, ..., z_{n-1})$.

The following lemma is a generalization of an approximation theorem in \cite{CDS1} and \cite{T2}.

\begin{lemma} \label{metrapp} Let $\omega\in c_1(L)$ be a  conical
K\"ahler current on $(X, \beta D)$ for some $\beta \in (0,1)$,
with bounded local potentials in $X$. Suppose that $\omega$  satisfies

\smallskip

{\rm (i)} $\omega^n =  e^{F} |s_D|_{h_D}^{-2\beta} \Omega$, where $\Omega
\in C^\infty(K_X\otimes \bar K_X)$, and $F\in C^0(X) \cap C^\infty(X\setminus D)$,

{\rm (ii)} $Ric(\omega) = -\ddbar \log \omega^n \geq -K \omega $ in the sense of currents, for some $K\in \mathbb{R}^+$.

\smallskip

Then $\omega\in C^\infty(X \setminus D)$ and there exist $K'>0$ and a sequence of smooth K\"ahler metrics $\omega_j \in c_1(L)$ such that

\smallskip

{\rm (1)} $Ric(\omega_j) \geq - K' \omega_j$ for all $j$;

{\rm (2)}
$\omega_j$ converges to $\omega$ smoothly on $X\setminus D$ and there exists $C>0$ such that for all $j$
$$\omega_j^n \leq C \omega^n;$$

{\rm (3)} Let $\omega_j = \omega+ \ddbar \varphi_j$ with $\sup_X \varphi_j=0$. Then $\varphi_j\to 0 $ uniformly in $L^\infty(X)$;

{\rm (4)} $(X, \omega_j)$ converges to the metric completion of $(X\setminus D, \omega)$ in Gromov-Hausdorff sense.

\end{lemma}

\medskip
\noindent
{\it Proof.}
Let $\hat \omega$ be a smooth K\"ahler metric in $c_1(L)$. Then there exists $\psi\in L^\infty(X)\cap PSH(X, \hat\omega)$ such that  $\omega= \hat\omega + \ddbar \psi$ due to the result of Kolodziej \cite{Kol}.
From the assumption,  we have
$$ - \ddbar F - \beta Ric(h_D) - \ddbar \log \Omega \geq - K (\hat \omega + \ddbar \psi)$$
in the sense of currents, and  the equation for $\psi$ is given by
\begin{equation}\label{eqnforpsi}
(\hat\omega+ \ddbar \psi)^n = e^{-(-F+ K\psi + \beta \log |s_D|_{h_D}^2) +K\psi} \Omega.
\end{equation}
Then $K\psi + \beta \log |s_D|^2_{h_D}- F \in PSH(X, K\hat\omega - \ddbar \log \Omega)$ and by Demailly's regularization theorem \cite{D}, there exist $A>0$ and a sequence of smooth
$$G_j\in PSH(X, A\hat\omega) \cap C^\infty(X)$$
 such that $G_j$ decreases to $K\psi + \beta \log |s_D|^2_{h_D}- F$.

Let $TX_z$ be the holomorphic tangent bundle of $X$ at $z$ and $exp\,h(z, \xi)$ be the quasi-holomorphic exponential map induced by the Chern connection of $\hat\omega$ constructed by Demailly \cite{D} for $\xi \in TX_z$.  The map $exp\,h_z(\xi)=exp\,h(z, \xi)$ is both smooth in a neighborhood of a fixed point $x\in X$ and $\xi\in \mathbb{C}^n$.   Let $\chi$ be cut-off function defined by the smooth function $\chi:\mathbb{R} \rightarrow \mathbb{R}$  satisfying
 $$\chi(t) >0, ~ t<1, ~\chi(t) =0, t\geq 1, ~\int_{\xi\in \mathbb{C}^n} \chi(|\xi|^2) d\lambda(\xi)=1. $$
 Then  $G_j$ can be explicitly expressed as in \cite{D} by
 \begin{eqnarray}
 \label{Demailly}
 G_j(z) &=& j^{2n} \int_{TX_z} (-F + K\psi+\beta \log |s_D|^2_{h_D})(exp\,h_z(\xi)) \chi(j^2 |\xi|^2) d\lambda (\xi) \label{appg}\\
  &=&   \int_{TX_z} (-F + K\psi+\beta \log |s_D|^2_{h_D})(exph_z(j^{-1}\xi)) \chi(|\xi|^2) d\lambda (\xi) , \nonumber%
 \end{eqnarray}
 where $d\lambda(\xi)$ is the standard Euclidean volume form on $\mathbb{C}^n$.
We observe that, for any point $z\in X\setminus D$, the support of the integrand in (\ref{Demailly}) is a sufficiently small neighborhood of $z$ for sufficiently large $j$, so if $\psi$ is in $C^k$ near $z$, then $G_j(z)$ will also be $C^k$ in a slightly smaller neighborhood of $z$. This formula makes sure that the support of the integrant does not include any point in $D$.

  We consider the following approximating Monge-Amp\`ere equations
\begin{equation}\label{appma}
(\hat\omega+\ddbar\psi_j)^n = e^{-G_j + K \psi_j} \Omega.
\end{equation}
Since $e^{-G_j}$ is uniformly bounded in $L^p(X, \Omega)$ for some $p>1$ as $\beta\in (0,1)$, the unique smooth solution $\psi_j$ is uniformly bounded in $L^\infty(X)$ \cite{Kol}.

We first derive a Ricci lower bound for $\omega_j=\hat\omega+\ddbar\psi_j$. There exists $A'>0$ such that for all $j$,
\begin{eqnarray}\label{lbd1}
Ric(\omega_j ) &=& -K(\hat\omega+\ddbar \psi_j) + K\hat\omega + \ddbar G_j - \ddbar\log \Omega \nonumber \\
&\geq& -K\omega_j - A' \hat\omega.
\end{eqnarray}
Then  there exists $C>0$ such that
\begin{equation}\label{ysch}
\omega_j= \hat\omega+ \ddbar \psi_j \geq C \hat\omega
\end{equation}
 for all $j$ by the argument of Yau's Schwarz lemma combined with the boundedness of $\psi_j$.  Combining (\ref{lbd1}) and (\ref{ysch}), there exists $K'>0$ such that for all $j$
\begin{equation}\label{lbd2}
Ric(\omega_j) \geq -K' \omega_j.
\end{equation}

 The diameter of $(X, g_j)$ is uniformly bounded above because there exists $C>0$ such that
\begin{equation}\label{c1alpha}
\omega_j = \hat\omega+ \ddbar \psi_j \leq C|s_D|^{-2\beta}_{h_D} \hat\omega
\end{equation}
for all $j$, after combining (\ref{ysch})  and  (\ref{appma}).  From (\ref{ysch}) and (\ref{c1alpha}), for any compact domain $\mathcal{K}\subset\subset X\setminus D$ and $\gamma \in (0, 1)$, there exists $C_{k, \mathcal{K}, \gamma}$ such that
\begin{equation}\label{c11}
\| \psi_j \|_{C^{1, \gamma}(\mathcal{K})} \leq C_{k, \mathcal{K}, \gamma}.
\end{equation}
Therefore by the uniqueness of the solution of (\ref{eqnforpsi}), $\psi_j$ converges to $\psi$.
By the formula (\ref{appg}) and the regularity of $\psi$ and $F$, for any compact domain $\mathcal{K} \subset\subset X\setminus D$, $G_j$ are uniformly bounded in $C^{1, \gamma}(\mathcal{K})$ for fixed $\gamma \in (0,1)$ for sufficiently large $j$ and so by Evans-Krylov's estimates,  $\psi_j$ is uniformly bounded in $C^{2, \gamma}(\mathcal{K})$. Therefore $\psi$ is uniformly bounded in $C^{2, \gamma'}(\mathcal{K})$ for $\gamma'\in (0, \gamma)$. By the formula (\ref{appg}) and linearizing (\ref{appma}), we can apply Schauder's estimates and bootstrap to derive uniform higher order estimates on $\psi_j$, $G_j$ and $\psi$ on $\mathcal{K}$. We then conclude that $\psi\in C^\infty(X\setminus D)$ and $\psi_j$ converges to $\psi$ on $C^\infty(X\setminus D)$.

Since $\psi_j$ converges smoothly on $X\setminus D$, we can fix a point $p\in X\setminus D$ such that the volume of the geodesic ball centered at $p$ of radius $1$ is uniformly bounded below away from $0$. Therefore we can apply Cheeger-Colding theory. Following the same arguments in \cite{Dat}, $(X, g_j)$ converges in Gromov-Hausdorff topology to the metric completion of $(X\setminus D, \omega)$ as $j\rightarrow \infty$. This completes the proof of the lemma.

\bigskip

We remark that Lemma \ref{metrapp} can be applied to the approximation of conical Fano K\"ahler-Einstein metrics in \cite{CDS1} and \cite{T2} without applying interpolation of the conical $F$-functionals. In particular, $X \setminus D$ is  convex and coincides with the regular part of the metric completion of $(X\setminus D, \omega)$.   Using this approximation lemma, we obtain the following Sobolev inequality:

\begin{lemma}
\label{Sobolev}
Let $(X, \beta D)$ and $\omega$ be defined as in the previous lemma. Then for any $f\in W^{1,2}(X)$,
$$\| f\|_{L^{n/(n-1)}(X, \omega) } \leq C(K, \ell) (\|\nabla f\|_{L^2(X, \omega)} + \|f\|_{L^2(X, \omega)} ),$$
where $\ell$ is the diameter of the metric completion of $(X\setminus D, \omega)$.

\end{lemma}

\noindent
{\it Proof.} Since $C^\infty(X)$ is dense in $W^{1,2}(X)$,
it suffices to prove the inequality for $f\in C^\infty(X)$.
Applying the approximation $\omega_j$ introduced in the previous lemma, we have the following Sobolev inequality for a uniform constant $C_S = C_S(K, \ell)>0$,
$$\| f\|_{L^{n/(n-1)}(X, \omega_j) } \leq C_S (\|\nabla f\|_{L^2(X, \omega_j)} + \|f\|_{L^2(X, \omega_j)} ),$$
for all $j$. Obviously as $j\rightarrow \infty$, we have
$$\| f\|_{L^{n/(n-1)}(X, \omega_j) } \rightarrow  \| f\|_{L^{n/(n-1}(X, \omega) }, ~  \| f\|_{L^{2 }(X, \omega_j) } \rightarrow \| f\|_{L^{2 }(X, \omega) } $$
as $\omega_j $ converges to $\omega$ uniformly in any $C^l(U)$ for any fixed $U\subset \subset X\setminus D$ and $\omega_j^n$ uniformly bounded above.  Suppose $\omega_j = \omega+ \ddbar \psi_j$ and $\|\psi_j \|_{L^\infty} \rightarrow 0$ as $j \rightarrow \infty$. Then
\begin{eqnarray}
\|\nabla f\|_{L^2(X, \omega_j)}^2 - \|\nabla f\|_{L^2(X, \omega)}^2
&=&  \int_X \partial f \wedge \dbar f \wedge \omega_j^{n-1} - \int_X \partial f \wedge \dbar f \wedge \omega^{n-1} \nonumber \\
&=& \int_X \partial f \wedge \dbar f \wedge (\omega_j - \omega) \wedge (\sum_{k=0}^{n-2} \omega^k \wedge \omega_j^{n-2-k}) \nonumber \\
&=&  \int_X \psi_j \ddbar f \wedge \ddbar f   \wedge (\sum_{k=0}^{n-2} \omega^k \wedge \omega_j^{n-2-k}).\nonumber
\end{eqnarray}
This last expression is clearly bounded in absolute value by $C_f\|\psi_j\|_{L^\infty}$ for some constant $C_f>0$. Letting $j\to\infty$, we obtain
$$\lim_{j\rightarrow \infty} \|\nabla f\|_{L^2(X, \omega_j)} =  \|\nabla f\|_{L^2(X, \omega)}  $$
and the lemma is proved.

\subsection{Monotonicity of $\mu(g,\tau)$}

First, we note that, under the conditions (\ref{Fspace}), the function $F^{1/2}$ is in $W^{1,2}(S^2,\beta)$, and has $L^2$ norm $1$. Thus the preceding Sobolev inequality (Lemma \ref{Sobolev}) holds for the function $F^{1/2}$, and applying as usual Jensen's inequality as in \cite{Rot} shows that
\begin{equation}
\mu(g,\tau)>-\infty.
\end{equation}
We claim next that $\mu(g(t),T-t)$ is increasing along the Ricci flow, just as in the case of smooth manifolds. For this we fix $T$, and consider the following coupled system of equations
\begin{eqnarray}
\ddt g(t)=-2\,Ric(g(t)),\qquad
{\partial F(t)\over\partial\tau}
=
\Delta_t F-R(t)F
\end{eqnarray}
with $\tau=T-t$. We need the following lemmas:

\begin{lemma}
\label{Fpreserve}
Let $F(t)$ evolve by the above coupled system.
If $F(t)$ satisfies the conditions (\ref{Fspace}) at $\tau=0$, then $F(t)$ will satisfy the conditions (\ref{Fspace}) for any $0\leq\tau\leq T$.
\end{lemma}

\noindent
{\it Proof}. First we note that the existence of a solution $F$ to the equation, satisfying the condition $F\in {\cal B}(S^2,\beta)$ for each $\tau$, is guaranteed by Yin's Theorem \ref{Yin2}, parts (1) and (2).  Recall that we have already proved in Section \S 2.2 that the class of regular metrics with conic singularities in the sense of Definition \ref{regularA} is preserved by the Ricci flow. In particular,
$\partial_t u$, $\partial_tR$ are in ${\cal B}((S^2,\beta)\times[0,T])$.
By part (3) of Theorem \ref{Yin2},
we also obtain $\partial_tF\in {\cal B}(S^2,\beta)$.
Substituting this into the backwards heat equation for $F(t)$, we find that $\Delta_tF(t)\in {\cal B}(S^2,\beta)$.

It remains to show that the condition $F\geq \delta_F>0$ is preserved along the flow. For this, it suffices to apply the maximum principle formulated in the next lemma. The proof of Lemma \ref{Fpreserve} is complete.

\bigskip
The next lemma is a maximum principle which holds thanks to the existence of a geometric barrier. Besides its use in the proof of Lemma \ref{Fpreserve}, it will also allow us to extend some classical arguments in the Ricci flow on smooth manifolds to the case of conic singularities.

\begin{lemma}
\label{maximum}
Let $g=g(t)$ be a $C^0$ metric on $(S^2\setminus \beta)\times [0,T]$.
Let $f\in C^2((S^2\setminus\beta)\times[0,T]) \cap L^\infty(S^2 \times[0,T])$ satisfy the following differential inequality
\begin{eqnarray}
\partial_t f\geq \Delta_g f+b(x,t) f
\end{eqnarray}
where $\Delta_g$ is the Laplacian with respect to $g$, and $b(x,t)$ is a bounded function.

{\rm (1)} If $f(x,0)\geq 0$ for all $x\in S^2\setminus\beta$, then $f(x,t)\geq 0$
in $(S^2\setminus\beta)\times[0,T]$;

{\rm (2)} If $f(x,0)\geq \delta>0$ for all $x\in S^2\setminus\beta$ and some constant $\delta>0$, then $f(x,t)\geq \delta e^{-t\|b\|_{L^\infty}}$
on $(S^2,\beta)\times[0,T]$. In particular, $f(x,t)>0$ on $(S^2\setminus\beta)\times[0,T]$.

\end{lemma}

\noindent
{\it Proof.} We prove (1) first. Replacing $f$ by $fe^{-At}$ for some large positive constant $A$, we can replace $b(x,t)$ by $b(x,t)-A$, and hence assume that $b$
is strictly less than $-B$, for any fixed positive constant $B$. Let $\tilde f=f-\epsilon\,\log\prod_{j=1}^k|s_j|_h^2$ where $s_j$ is a holomorphic section of $K_{S^2}^{-1}$
with divisor $[2p_j]$, normalized so that $|s_j|_h^2\leq {1\over 2}$, with $h$ a smooth metric on $K_{S^2}^{-1}$. The function $\tilde f$ satisfies the differential inequality
\begin{eqnarray}
\partial_t\tilde f\geq \Delta_g \tilde f+b\tilde f +\epsilon(k\,R_h+b\log\prod_{j=1}^k|s_j|_h^2)
\end{eqnarray}
and thus, since $R_h$ is the contraction of the curvature of $h$ with $g$ and hence bounded,
\begin{eqnarray}
\partial_t\tilde f
\geq
\Delta_g \tilde f+b\tilde f,
\end{eqnarray}
if we choose $B$ to be sufficiently large. Since $\tilde f\to +\infty$ near each of the conical singularities $p_j$, it must attain its minimum somewhere in $(S^2\setminus\beta)\times [0,T]$. Assume that this minimum is strictly negative, and let $t_0>0$ be the first time when it is achieved, at some point $x_0\in S^2\setminus\beta$. The above differential inequality would imply that
\begin{eqnarray}
\partial_t\tilde f(x_0,t_0)
\geq b(x_0,t_0)\tilde f(x_0,t_0)>0.
\end{eqnarray}
But this would imply in turn that $\tilde f$ must have attained values strictly lower than $\tilde f(x_0,t_0)$, which is a contradiction. Thus the minimum of $\tilde f$ must be non-negative. Letting $\epsilon\to 0$, it follows that $f$ is non-negative, and (1) is proved.

\smallskip
Next, we prove (2). This time, we set $\tilde f=f-\epsilon e^{-At}$. Then the function $\tilde f$ satisfies the differential inequality
\begin{eqnarray}
\partial_t\tilde f\geq \Delta_g\tilde f+b\tilde f
+\epsilon(A+b)e^{-At}
\geq
\Delta_g\tilde f+b\tilde f
\end{eqnarray}
for $A=\|b\|_{L^\infty}$. In view of Part (1), we have $\tilde f\geq 0$ for
all $t\in [0,T]$, if $\tilde f\geq 0$ at $t=0$.  Thus we choose $\epsilon=\delta$,
and obtain the bound $f(x,t)\geq \delta e^{-t\|b\|_{L^\infty}}$,
as claimed. The proof of Lemma \ref{maximum} is complete.

\bigskip

We can complete now the proof of the monotonicity of the function $\mu(g,\tau)$ along the Ricci flow exactly as in Perelman's original arguments (c.f. \cite{Top}). Thus set
\begin{equation}
\label{vv}
v=\bigg(\tau (2\Delta f-|\nabla f|^2+R)+f-2\bigg){e^{-f}\over 4\pi\tau}.
\end{equation}
Note that $v$ differs from the integrand in the functional $W(g,f,\tau)$ by
$2\Delta F$. Since $F$ is bounded and has finite Dirichlet energy, we have
$\int_{S^2}\Delta F=0$, and thus
\begin{equation}
W(g,f,\tau)
=
\int_{S^2} vdV_t.
\end{equation}
Now the same calculation on the smooth part $S^2\setminus\beta$ as in Perelman's original arguments gives
\begin{equation}
\Box^* v = - 2\tau |Ric(g) + Hess(f) - (2\tau)^{-1} g |^2 {e^{-f}\over 4\pi\tau},
\end{equation}
where $\Box^*= -\ddt{} - \Delta_t+R(t)$. It follows that
\begin{eqnarray}
\label{Wv}
\ddt {W}(g, f, \tau)
=\ddt{}\int_{S^2} v dV_t
=\int_{S^2} ( -\Box^* v- \Delta_t v) dV_t .
\end{eqnarray}

\begin{lemma}
\label{v}
The function $v$ is bounded and has bounded Dirichlet energy.
\end{lemma}

\noindent
{\it Proof}. Indeed, $v$
can be rewritten in terms of $F={e^{-f}\over 4\pi\tau}$ as
\begin{equation}
\label{v1}
v=-2\tau \Delta F+\tau{|\nabla F|^2\over F}
+
(\tau R+f-2)F.
\end{equation}
Since $F\geq \delta_F>0$, it is clearly that the last term on the right is bounded and has bounded Dirichlet energy. Clearly $\Delta F$ is bounded and has bounded Dirichlet energy since $\Delta F\in {\cal B}(S^2,\beta)$. Since $F\in {\cal B}(S^2,\beta)$ and $\Delta F$ is bounded, the fact that $|\nabla F|^2$ is bounded and has bounded energy follows from Lemma 3.2 in \cite{Yi2}. Lemma \ref{v} is proved.

\bigskip
As a consequence of Lemma \ref{v}, we have
\begin{equation}
\int_{S^2}\Delta_t v\,dV_t=0,
\end{equation}
and hence
\begin{eqnarray}
 \ddt {W}(g, f, \tau)
= \int_{S^2}  |Ric(g) + Hess(f)  - (2\tau)^{-1} g|^2 {e^{-f} \over 4\pi\tau}
dV_t\geq 0.
\end{eqnarray}
Similarly, the normalization (\ref{Fnormalization}) is preserved along the coupled equations.

\bigskip

We have now all the ingredients for the proof of the monotonicity of $\mu(g(t),\tau)$. Let $t_1<t_2$. For any function $F\in {\cal B}(S^2,\beta)$,
$F\geq \delta_F>0$, $\Delta_{t_2}F\in {\cal B}(S^2,\beta)$,
let $F(t_1)$ be the solution at time $t_1$ of the backwards heat equation, with final data $F(t_2)=F$. By Lemma \ref{Fpreserve}, $F(t_1)$ satisfies all the conditions
(\ref{Fnormalization}). Then by the monotonicity
of the $W$-functional along the Ricci flow coupled with the backwards heat equation,
\begin{eqnarray}
W(g(t_2),F,\tau(t_2))
=
W(g(t_2),f(t_2),\tau(t_2))
\geq
W(g(t_1),f(t_1),\tau(t_1))
\geq
\mu(g(t_1),\tau(t_1)).
\end{eqnarray}
Taking the infimum over all $F$ satisfying the conditions (\ref{Fspace}) and (\ref{Fnormalization}) gives the desired monotonicity
\begin{eqnarray}
\mu(g(t_2),\tau(t_2))
\geq
\mu(g(t_1),\tau(t_1)),
\qquad t_1\leq t_2.
\end{eqnarray}

\subsection{Estimates for the curvature $R(t)$}

Once we have the monotonicity of $\mu(g,\tau)$ along the Ricci flow,
we can apply Perelman's arguments in the same way \cite{SeT}
and obtain the following:

\begin{lemma} \label{pere}
Let $g(t)$ be the Ricci flow with an initial metric $g_0$ which is a regular metric with conical singularities in the sense of Definition \ref{regularA}.
Let $ v(t) $ be the Ricci potential of $g(t)$ defined by
\begin{equation}
{\bf Ric}(g(t)) =
{1\over 2}\chi(S^2,\beta)
g(t)  - \ddbar v(t) +  \sum_{i=1}^k \beta_i [p_i]
\end{equation}
and $ \int_{S^2}e^{-v (t)} dg(t) = 2.$ Then there exists $C=C(g_0)>0$ such that for all $t\geq 0$ and $p\in S^2\setminus\beta$
\begin{equation}
\label{perel} (|v| + |\nabla v |_{g} + |\Delta v|_g)(p, t) \leq C,
\end{equation}
In particular, the curvature of $g$ is uniformly bounded along the flow.
Furthermore, the diameter of $S^2\setminus\beta$ with respect $g(t)$
is uniformly bounded, and the flow satisfies Perelman's $\kappa$-noncollapsing condition.

\end{lemma}

\noindent
{\it Proof}. Since the Ricci flow preserves regularity in the sense of Definition \ref{regularA}, and since the initial metric $g_0$ is regular, it follows that $R(t)$ is a bounded function on $S^2\setminus\beta$ for each $t$. This implies that
$\Delta v$ is a bounded function for each $t$. By Theorem 3.1 of \cite{Yi1}, the function $v$ is also bounded and has bounded Dirichlet energy.

All the arguments of Perelman in \cite{SeT} can be adapted to the conical case
using the monotonicity of the modified conical $W$-functional.
All we need to do is to justify how to apply the maximum principle to derive the uniform bounds on $v$. By the same argument in \cite{SeT}, $v$ is uniformly bounded below and we can assume that $v+B\geq 1$ for all $t\in [0, \infty)$. Let $H= \frac{ |\nabla v|^2 }{v+B}$. Then the calculations in \cite{SeT} show that
$$ \ddt{H}\leq \Delta H + |\nabla H|^2+ A (v+B)^{-1} H - \delta (v+B)^{-1} H^2$$
for some $A, \delta > 0$
We cannot directly apply the maximum principle to obtain a uniform upper bound for $H$. As in Lemma \ref{maximum}, we let $\sigma$ be a holomorphic section in $-kK_{S^2}$ such that $\sigma$ vanishes at each $p_i$ to order $2$. We consider $f_{\epsilon}$ defined by
$$ f_\epsilon = e^H + \epsilon \log |\sigma|^2_{\omega^k} - \epsilon k t,$$
where $\omega$ is the volume form of $g(t)$.
Straightforward calculations show that
$$ \ddt {f_\epsilon} \leq \Delta f_\epsilon + \frac{Ae^H}{v+B} H - \frac{\delta e^H}{v+B} H^2$$
since $(\ddt{} - \Delta) \log |\sigma|^2_{\omega^k} = k(1-{1\over 2}\sum_{i=1}^k\beta_i)$.
We can assume $\sup_{S^2\times [0, T]} f_\epsilon = f_\epsilon(q, t)$ for fixed $T>0$ since $f_\epsilon$ tends to $-\infty$ near each conical point. Then $H(q, t) \leq {A/\delta}$ by the maximum principle, since $p$ cannot be a conical point. Immediately we have for any $(q', t')\in (S^2\setminus \beta)\times [0, T]$,
$$e^{H(q', t')} \leq e^{A/\delta}  + \epsilon \sup_{S^2\times [0, T] } \log |\sigma|^2_{\omega^k}- \epsilon \log |\sigma|^2_{\omega^k} (q', t') + \epsilon k T. $$
By letting $\epsilon\rightarrow 0$, we obtain
\begin{equation}\label{grades}
H \leq  \frac{A}{ \delta}
\end{equation} on $S^2 \times [0, \infty).$
A similar argument shows that $\frac{\Delta u}{v+B}$ is also uniformly bounded.
The estimate (\ref{perel}) is proved.

The same argument of Perelman \cite{SeT} can be applied to obtain Perelman's $\kappa$-nonlocal collapsing and then a uniform diameter bound for $(S^2, g(t))$ for all $t>0$ due to the estimate (\ref{grades}) and the geodesic convexity of $(S^2\setminus \beta)$ in $(S^2, g(t))$. The proof of the lemma is complete.

\bigskip

One does not expect the global Hamilton-Cheeger-Gromov compactness to hold in this case since the injectivity radius will in general fail to be bounded away from $0$. In fact, in the semi-stable and unstable case, for generic points, the injectivity radius will always tend to $0$ as the gauge transformations are the $\mathbb{C}^*$ action, moving all but $p_k$ to the other limiting conical point. Such phenomena will be addressed in the following sections when Troyanov's stability condition is violated.


\section{
The stable case }

In this section we prove the convergence for the conical Ricci flow on $(S^2, \beta)$ in the stable case. The following is a more precise version of Part 1, Theorem \ref{main}:

\begin{lemma} If $2\beta_{max} < \sum_{i=1}^k   \beta_i$, then for any regular initial metric $g_0=e^{u_0} g_\beta \in c_1(S^2)$,
the conical Ricci flow converges to the unique constant curvature metric $g_\infty \in c_1(S^2)$ on $(S^2, \beta)$ in the sense that the potentials $\varphi$
are uniformly bounded in some Schauder space
$C^\alpha(S^2)$ for some $\alpha>0$,
and they converge in $C^\infty$ on any compact subset of $S^2\setminus\beta$. Furthermore, $(S^2, g(t))$ converges in Gromov-Hausdorff topology to the unique constant curvature metric $g_\infty$ on $(S^2, \beta)$.

\end{lemma}

\noindent
{\it Proof.} The proof is an adaptation of the methods in \cite{PS} and \cite{PSS}, exploiting the properness of the functional $F_\beta$ and the fact, established in the previous section, that the Ricci potential $v$ is uniformly bounded along the K\"ahler-Ricci flow.

First, we note as in \cite{PSS} that if we express the metrics $g(t)$ along the K\"ahler-Ricci flow as
\begin{equation}
g(t)=g_{FS}+{\sqrt{-1}\over 2\pi}\partial\bar\partial\varphi(t)
\end{equation}
then $\dot\varphi(t)+v=c(t)$, where $c(t)$ is a constant depending only on the time $t$. Arguing as in \cite{PSS}, we see that the constant $c(t)$ can be made uniformly bounded in $t$ by choosing suitably the arbitrary constant in the definition of $\varphi(0)$. The integration by parts in the argument are justified because $v$ is bounded and has bounded Dirichlet energy.
With this choice of normalization, we have then
\begin{equation}
{\rm sup}_t\|\dot\varphi\|_{C^0} <\infty.
\end{equation}
This estimate together with the properness of $F_\beta$ can now be shown to imply the following key estimate for the gradient
\begin{equation}
\label{gradient}
\int_{S^2}{\sqrt{-1}\over 2\pi}
\partial\varphi\wedge\bar\partial\varphi \leq C
\end{equation}
and for the average of $\varphi$ along the flow
\begin{equation}
\label{average}
|\int_{S^2}\varphi\,dg_{FS}|\leq C.
\end{equation}
To see this, we begin by noting, as in \cite{PS}, that the functional $F_\beta$ and the functional $F_\beta^0$ defined by
\begin{equation}
F_\beta^0(\varphi)={i\over 8\pi}\int_{S^2}\partial\varphi\wedge\bar\partial\varphi
-{1\over 2}\int_{S^2}\varphi\,dg_{FS}
\end{equation}
are comparable along the flow,
\begin{equation}
|F_\beta(\varphi)-F_\beta^0(\varphi)|\leq C.
\end{equation}
This is because their difference satisfies
\begin{equation}
2e^{-\|\dot\varphi\|_{C^0}}
\leq |\int_{S^2}e^{-{1\over 2}\chi(S^2,\beta)\varphi}
\prod_{i=1}^k
|\sigma_i|_{\omega_{FS}}^{-\beta_i}
g_{FS}
|
=
|\int_{S^2}e^{-\dot\varphi}(\omega_{\beta}+{\sqrt{-1}\over 2\pi}\partial\bar\partial\varphi)|
\leq 2e^{\|\dot\varphi\|_{C^0}}
\end{equation}
which is uniformly bounded since $\|\dot\varphi\|_{C^0}$ is uniformly bounded.

Next, a straightforward calculation shows that $F_\beta$ is decreasing along the K\"ahler-Ricci flow, and hence, using the properness of $F_\beta$,
\begin{equation}
F_\beta(\varphi(0))
\geq F_\beta(\varphi)\geq C_1\int_{S^2}{\sqrt{-1}\over 2\pi}\partial\varphi\wedge\bar\partial\varphi-C_2
\end{equation}
which shows that
\begin{equation}
\int_{S^2}{\sqrt{-1}\over 2\pi}\partial\varphi\wedge\bar\partial\varphi
\leq C_3
\end{equation}
for all $t$, which is equation (\ref{gradient}). It follows also that $|F_\beta(\varphi)|$ is uniformyly bounded, and hence that $|F_\beta^0(\varphi)|$ is uniformly bounded.
Since we already know that (\ref{gradient}) holds, the estimate (\ref{average}) follows at once.

\smallskip
We can now apply the Trudinger inequality on compact Riemannian manifolds of dimension $2$: there exists constants $C>0$ and $\kappa>0$ so that for any $p>0$,
\begin{equation}
\label{Trudinger}
\int_{S^2}
e^{p|\varphi|}dg_{FS}
\leq
C
{\rm exp}\bigg(\kappa p^2\int_{S^2}{\sqrt{-1}\over 2\pi}\partial\varphi\wedge\bar\partial\varphi+p|({1\over 2}\int_{S^2}\varphi\,dg_{FS}|\bigg)
\end{equation}
We deduce that for any $p$,
\begin{equation}
{\sup}_t\int_{S^2}e^{p|\varphi|}dg_{FS}
<\infty.
\end{equation}
Next we rewrite the equation for the flow as
\begin{equation}
g_{FS}
+
{\sqrt{-1}\over 2\pi}
\partial\bar\partial\varphi
=
g_{FS} e^{\dot\varphi}
e^{-{1\over 2}\chi(S^2,\beta)\varphi}
\prod_{i=1}^k
|\sigma_i|_{g_{FS}}^{-\beta_i}
\end{equation}
Since $e^{{1\over 2}\chi(S^2,\beta)|\varphi|}$ is in $L^p$ for any $p<\infty$, we can apply H\"older's inequality and find that the right hand side of the above equation is in $L^p$ for some $p>1$. By the standard $W^{2,p}$ estimate for elliptic linear PDE, we can conclude that
\begin{equation}
\varphi-{1\over 2}\int_{S^2}\varphi dg_{FS}
\end{equation}
is uniformly bounded in $W^{2,p}(S^2)\subset C^\alpha(S^2)$
for $\alpha=2-{2\over p}$. In view of (\ref{average}), we can conclude that
$\varphi$ is uniformly bounded in $C^\alpha(S^2)$.

\smallskip
We can then apply the parabolic versions (\cite{PSS}) of the
standard estimates of Aubin and Yau for the second derivatives, and of Calabi for the third order derivatives, modified by $\epsilon \log
\prod_{i=1}^k|\sigma_i|_{FS}^2$, to obtain the uniform boundedness of the potentials $\varphi$
in $C^\infty(K)$ for any compact subset $K\subset\subset S^2\setminus\beta$.

\smallskip
Applying the arguments in \cite{PSSW} (c.f. Lemma \ref{curcon}), one can show that $\| \dot \varphi(t) \|_{L^\infty} \rightarrow 0$ as $t\rightarrow \infty$. This implies the convergence of the potentials $\varphi$ in $C^\infty$ on compact subsets of $S^2\setminus\beta$ for a subsequence $\varphi(t_j)$. Since the limit is unique, the convergence along subsequences implies the convergence of the whole flow. The Gromov-Hausdorff convergence follows immediately because $\|\dot\varphi\|_{C^\alpha(S^2, g_{FS})}$ and $\|\varphi\|_{C^\alpha(S^2, g_{FS})}$ are uniformly bounded for all $t>0$ for some fixed $\alpha>0$ and the standard sphere metric $g_{FS}$. The proof of the lemma is complete.


\section{The semi-stable and unstable cases}

In this section, we give the proof of Theorem \ref{main} in the semi-stable and stable cases.

\subsection{Sequential convergence}

The sequential convergence for the conical Ricci flow in both the semi-stable and unstable cases can be established as follows.

\begin{lemma}\label{ghcon} After passing to a subsequence, for any sequence $ t_j \rightarrow \infty$, the spaces
$(S^2, g(t_j))$ converges to a compact length metric space $(X, d)$ satisfying the following properties:

{\rm (1)} $X$ is homeomorphic to $S^2$;

{\rm (2)} the singular set $D$ is a finite set of isolated points;

{\rm (3)} $(X\setminus D, d)$ is a smooth surface equipped with a smooth Riemannian metric with volume $2$.

In particular, the convergence is smooth on $X\setminus D$.
\end{lemma}

\noindent
{\it Proof}.   First we can apply Lemma \ref{metrapp} and obtain a constant $K>0$ and smooth metrics $g_j$ for all $j$ such that
$$R(g_j)\geq - K, ~ d_{GH}( (S^2, g(t_j)), (S^2, g_j)) \leq  j^{-1} $$
In particular, $g_j\in c_1(S^2)$ and the diameter of $(S^2, g_j)$ is uniformly bounded. We can now directly apply Cheeger-Colding theory and obtain a Gromov-Hausdorff limit $(X, d)$, a compact metric length space, after taking a convergent subsequence of $(S^2, g_j)$.  Without loss of generality and by passing to a convergent subsequence, we can assume that $(S^2, g(t_j))$ converges to $(X, d)$ in Gromov-Hausdorff topology.

\smallskip
 We would like to show that the singular set $D$ of $(X, d)$ is finite and that in fact, it coincides with the set of limiting points of  the conical points along the sequence $(S^2, g(t_j))$. We let $D'\subset D$ be the set of  all the limiting points of the conical points. Obviously, $D'$ must be finite. We will show that $D=D'$. Suppose $P\in X\setminus D'$,  then there exist a sequence of points $P_j \in (S^2, g(t_j))$ converging to $P$ in Gromov-Hausdorff sense. Since $D'$ is finite, we can assume that the distance from $P$ to $D'$ is  bounded from below by $2r>0$. Therefore the distance from $P_j$ to the set of all conical points in $(S^2, g(t_j))$ is bounded from below by $r$ for sufficiently large $j$.  We then consider the sequence of balls $B_{g(t_j)}(P_j, r)$, which do not contain any conical point. Since the curvature of $g(t_j)$ is uniformly bounded on $B_{g(t_j)}(P_j, r)$ and one has uniformly nonlocal $\kappa$-collapsing for all $(S^2, g(t_j))$, the injectivity radius of any point in $B_{g(t_j)}(P_j, r/A)$ is uniformly bounded below by applying Klingenberg's lemma (see section 8.4 in \cite{Top}) for a fixed sufficiently large $A>0$.  Therefore, $B_{g(t_j)}(P_j, r/A)$ converges in $C^{1,\alpha}$ after passing to a convergent subsequence by Anderson's results. This implies that $P$ must be a regular point of $X$ and so $D=D'$.

\smallskip
 We can now establish the partial $C^0$-estimates as in \cite{DS}. Of course, one has to make a slight modification and this is essentially the case  studied in \cite{CDS2, T2}. However, in our situation, the singular set is much simpler since there is no singular set of Hausdorff codimension greater  than $2$ and thus each tangent cone of $(X, d)$ is a flat metric cone on $\mathbb{C}$. This implies that each tangent cone  is good, i.e., one can construct appropriate cut-off functions. Hence one can immediately obtain the partial $C^0$-estimates for the evolving metric $g(t)$. More precisely, there exists $\epsilon>0$ and $N>0$ such that for all $t\geq 0$ and any $p\in S^2$, there exists $\sigma \in H^0(S^2, K_{S^2}^{-N})$ satisfying
$$\left( |\sigma|^2 (\omega(t))^N \right) (p) \geq \epsilon, ~ \int_{S^2} |\sigma|^2 (\omega(t))^{N+1}=2.$$
Here $\omega(t)$ is the volume form of $g(t)$ and so it is a hermitian metric on $K_{S^2}^{-1}$. Suppose $p_\infty$ is a singular point in $(X, d)$. Then any tangent cone at $p_\infty$ must be a metric cone $\mathbb{C}_\gamma$ on $\mathbb{C}$ with a cone metric $g_\gamma=\ddbar |z|^{2\gamma}$ for some $\gamma \in (0, 1]$ with $0$ being $p_\infty$. The trivial line bundle on $\mathbb{C}_\gamma$ is equipped with the hermitian metric $e^{- |z|^{2\gamma}}$. Let $0\leq F\leq 1$ be the standard smooth cut-off function on $[0, \infty)$ with $F=1$ on $[0, 1/2]$ and $F=0$ on $[1, \infty)$. We then let
$$ \rho_\epsilon =F\left({\eta_\epsilon\over \log \epsilon} \right), ~ \eta_\epsilon = \max(\log |z|^2, 2\log \epsilon).$$
Then one can show by straightforward calculations  that
$$\int_{\mathbb{C}}|\nabla \rho_\epsilon|^2 dg_\gamma < C(- \log \epsilon)^{-1}$$
for some fixed constant $C>0$ uniform in $\epsilon\in (0, 1]$.
Obviously for any $K \subset\subset \mathbb{C}^*$ and $\delta>0$, there exists sufficiently small $\epsilon>0$ such that
$ \rho_\epsilon =1$ on $K$, ${\rm supp}\, \rho_\epsilon \subset\subset \mathbb{C}^*$.
Using the construction of $\rho_\epsilon$, one can prove the partial $C^0$-estimate as in \cite{DS}.

\smallskip

We can now make use of the arguments of \cite{DS}: the existence of the above sections $\sigma(\omega(t))$ for any $t$ implies that the surfaces $(S^2,\omega(t))$ can be uniformly imbedded into some ${\bf CP}^N$, separating points, and that the limit of their images must be a normal variety.
Since this normal variety is a projective degeneration of $S^2$, it must be $S^2$.

\smallskip

From Shi's local estimates, namely that uniform bounds for the curvature along the Ricci flow on any bounded domains implies similar bounds for the derivatives of the curvature on smaller domains,
the convergence on $X\setminus D$ is smooth and the limiting metric is a smooth metric on $X\setminus D$. The lemma is then proved.

\bigskip

Xiaochun Rong has pointed out to us that that one can apply Perelman's stability theorem for Alexandrov spaces instead of the partial $C^0$-estimate to show that $X$ is homeomorphic to $S^2$. But what was established above via
the partial $C^0$ estimate is slightly stronger: the image in ${\bf P}^N$ is a smooth ${\bf P}^1$, excluding the possibility of singular rational curve which is also homeomorphic to $S^2$, e.g. a rational curve with cuspidal singularities.

\bigskip
To identify the metric on the limiting space, we make use next of Hamilton's entropy functional \cite{H}, defined for metrics with ${\rm inf}_{S^2}R>0$ by
\begin{equation}
N= \int_{S^2} R \log R\, dg.
\end{equation}
 The assumption $\inf_{S^2} R_0>0$ can be removed by a trick of Chow \cite{Ch},
 which still works in exactly the same way in the case of the sphere with marked points (c.f. section 8.2, Chapter 5 \cite{CK}), after  replacing $R$ by $ R - s$, where $ s$ is defined by $\ddt{s} = s(s- {1\over 2}\chi(S^2,\beta))$ with $s(0) < \inf_{S^2} R_0$. Thus we can henceforth assume that ${\rm inf}_{S^2}R>0$.

\begin{lemma}
\label{entropy}
Let $v$ be the Ricci potential defined as in Lemma \ref{pere}.  Along the conical Ricci flow, we have
\begin{equation}
\ddt{N} = - \int_{S^2} \frac{|\nabla R + R \nabla v|^2}{R} d g - 2 \int_{S^2} |\nabla\nabla v - \frac{1}{2} (\Delta v) g |^2 d g,
\end{equation}
if $\inf_{S^2}R_0>0$.
\end{lemma}

\noindent
{\it Proof}. From the flow equation for $R$ and the maximum principle Lemma \ref{maximum}, it follows that $\inf_{S^2}R>0$ for all $t$ if $\inf_{S^2}R_0>0$. Thus the entropy functional $N$ is well-defined for all time. It suffices now to apply the same arguments as in Hamilton \cite{H}. The integration by parts which are required are justified in the lemmas which we state below.

\begin{lemma} Let $g=e^ug_\beta$ be any regular metric
in the sense of Definition \ref{regularA}.
Then
$$\int_{S^2} |\nabla^2 u|^2 dg= \int_{S^2} (\Delta u)^2 dg+ \int_{S^2} R^i_j \nabla^i u\nabla_j u ~dg,$$
where $\Delta$ is the Laplacian with respect to $g$.
In particular, $\int_{S^2} |\nabla^2 u|^2 dg < \infty $.

\end{lemma}

\noindent
{\it Proof.} First we notice that $\Delta u = - R + e^{-u} R_\beta$ and so it is bounded. By Yin's estimate for the Poisson equation, $|\nabla u|$ is bounded
because $u$ has bounded Dirichlet energy.
We now pick a family of cut-off functions  $\rho_\epsilon$,
$0\leq \rho_\epsilon\leq 1$, with the following properties: for any $\epsilon>0$ and any $K\subset\subset S^2\setminus\beta$, $\rho_\epsilon\in C_0^\infty(S^2\setminus \beta)$, with $\rho_\epsilon=1$ on $K$, and
$\int_{S^2} |\nabla\rho_\epsilon|^2\omega_\beta<\epsilon$.
\begin{eqnarray*}
\int_{S^2} \rho_\epsilon^2 \nabla^i \nabla^j u \nabla_i \nabla_j u
&=& -\int_{S^2} \rho_\epsilon \nabla^i \rho_\epsilon \nabla^j u \nabla_i \nabla_j u - \int_{S^2} \rho_\epsilon^2~  \nabla^j u \nabla^i \nabla_j \nabla_i u\\
&=&  -\int_{S^2} \rho_\epsilon \nabla^i \rho_\epsilon \nabla^j u \nabla_i \nabla_j u - \int_{S^2}\rho_\epsilon^2  \nabla^j u \nabla_j \nabla^i \nabla_i u + \int_{S^2} \rho_\epsilon^2 R^i_j \nabla^j u \nabla_i u\\
&=& \int_{S^2} \rho_\epsilon^2 (\Delta u)^2+ \int_{S^2} \rho_\epsilon^2 R^i_j \nabla^j u \nabla_i u - \int_{S^2} \rho_\epsilon \nabla^i \rho_\epsilon \nabla^j u \nabla_i \nabla_j u \\
&&
\qquad  + \int_{S^2} \rho_\epsilon \nabla_j \rho_\epsilon \nabla^j u  \Delta u.
\end{eqnarray*}
Hence
\begin{eqnarray*}
&& \left| \int_{S^2} \rho_\epsilon^2 |\nabla^2 u|^2  - \int_{S^2} \rho_\epsilon^2 (\Delta u)^2- \int_{S^2} \rho_\epsilon^2 R^i_j \nabla^j u \nabla_i u              \right |\\
&\leq & \left|  \int_{S^2} \rho_\epsilon \nabla^i \rho_\epsilon \nabla^j u \nabla_i \nabla_j u\right| + \left|     \int_{S^2} \rho_\epsilon \nabla_j \rho_\epsilon \nabla^j u  \Delta u)\right|   \\
&\leq& \left( \int_{S^2} |\nabla u|^2 |\nabla \rho_\epsilon|^2 \right)^{1/2} \left( \int_{S^2} \rho_\epsilon^2 |\nabla^2 u|^2 \right)^{1/2} + \left( \int_{S^2} |\nabla \rho_\epsilon|^2\right)^{1/2} \left( \int_{S^2} |\nabla u||\nabla u|\right)^{1/2}\\
&\leq &\left( \int_{S^2} |\nabla u|^2 |\nabla \rho_\epsilon|^2 \right)^{1/2}  \left(\int_{S^2} \rho_\epsilon^2 |\nabla^2 u|^2  - \int_{S^2} \rho_\epsilon^2 (\Delta u)^2- \int_{S^2} \rho_\epsilon^2 R^i_j \nabla^j u \nabla_i u      \right)^{1/2} \\
&&+ C\left( \int_{S^2} |\nabla u|^2 |\nabla \rho_\epsilon|^2 \right)^{1/2} + \left( \int_{S^2} |\nabla \rho_\epsilon|^2\right)^{1/2} \left( \int_{S^2} |\nabla u||\nabla u|\right)^{1/2}
\end{eqnarray*}
for some uniform constant $C>0$ since $|\nabla u|$, $\Delta u$ and $R$ are bounded. Then the proposition is proved by letting $\epsilon \rightarrow 0$.

\medskip

More generally, we have the following proposition, which can be proved in exactly the same way:

\begin{lemma} Let $g$ be a regular metric  $g=e^u g_\beta$ in the sense of Definition \ref{regularA}.
Suppose $f\in {\cal B}(S^2,\beta)$ with $\Delta f$ having bounded Dirichlet energy.
Then
$$\int_{S^2} |\nabla^2 f |^2 dg= \int_{S^2} (\Delta f)^2 dg+ \int_{S^2} R^i_j \nabla^i f \nabla_j f ~dg,$$
where $\Delta$ is the Laplacian with respect to $g$.
In particular, $\int_{S^2} |\nabla^2 f|^2 dg < \infty $.

\end{lemma}

Since $R\log R$ is bounded from below, the entropy $N$ is bounded from below, and Lemma \ref{entropy} implies immediately

 \begin{lemma}\label{entropy1} If $\inf_{S^2} R_0>0$, then

 \begin{equation}
 \lim_{t\rightarrow \infty} \int_{s=t}^{t+1}\int_{S^2} |\nabla\nabla v - \frac{1}{2} (\Delta v) g |^2 d g ~ds =0
 \end{equation}

 \end{lemma}

\medskip

The following lemma is  well-known in complex analysis and it will help establish the limiting soliton equation of the conical Ricci flow.

\begin{lemma}\label{onevar} Let $f(z, \bar z)$ be a smooth real-valued harmonic function on the punctured disc $\mathbb{B}^*\subset \mathbb{C}$. Then
$$f(z, \bar z) = Re(F(z)) + c \log |z|^2,$$
where $F(z)$ is a holomorphic function on $\mathbb{B}^*$ and $c\in \mathbb{R}$. In particular if $e^f \in L^1(\mathbb{B})$, then $F$ extends to a holomorphic function on $\mathbb{B}$.

\end{lemma}

\noindent
{\it Proof}. Let $h$ be the exponential map from the left plane $\{ w\in \mathbb{C}~|~Re (w) <0\}$ to $\mathbb{B}^*$. Then $u(w)= h^*f(w)= f( e^w)$ is also a harmonic function satisfying $u(w+2\pi \sqrt{-1})= u(w)$. Since the left plane is simply connected, there exists a complex conjugate $v(w)$ for $u(w)$. In particular,  $\nabla v(w+2\pi \sqrt{-1}) = \nabla v(w)$ for all $w$ and we can define a holomorphic function $G$ satisfying
$$G(w) = u(w) + \sqrt{-1} v(w) - cw, ~ G(w+2\pi i\sqrt{-1}) = G(w)$$
 for some $c\in \mathbb{R}$. Hence
 $$ f(z, \bar z) = Re ( G(\log z) )+ c \log |z|, $$
and $F(z)= G(\log z)$ is obviously a holomorphic function on $\mathbb{B}^*$ because $G(w+2\pi \sqrt{-1})= G(w)$.

If $e^f\in L^1(\mathbb{B})$,  then
$$\int_{\mathbb{B}} e^f = \int_{\mathbb{B}} |z|^{2c} |e^{F/2}|^2 <\infty.$$ Hence $z^m e^{F/2}$ is a holomorphic function on $\mathbb{B}$ for some sufficiently large $m\in \mathbb{Z}^+$ and this implies that $F$ cannot have a singularity at $0$. The proof of the lemma is complete.

\begin{lemma}
 \label{L1}
 Let $t_j \rightarrow \infty$. Then by passing to a subsequence,
 $(S^2, g(t_j))$ converges in Gromov-Hausdorff topology to one of the following:
 \begin{enumerate}

 \item a conical metric space $(S^2, \beta_\infty)$ of constant curvature $1-\frac{1}{2} \sum_{i=1}^k \beta_i$,

 \item
 a rotationally symmetric conical shrinking gradient Ricci soliton on $(S^2, \beta_\infty)$ with  $\beta_\infty= \beta_{p_\infty} [p_\infty] + \beta_{q_\infty} [q_\infty]$ with $0\leq \beta_{q_\infty} < \beta_{p_\infty} <1$.

\end{enumerate}

\end{lemma}

\noindent
{\it Proof}.
Using Lemma \ref{pere}, Lemma \ref{ghcon}  and Shi's local curvature estimates, we see that the flow converges smoothly on  $X\setminus D$. Since $R = \Delta v + (1-\frac{1}{2} \sum_{i=1}^k \beta_i)$ and $v$ is uniformly bounded in $C^0$, it
follows from the standard  estimates for the Laplace equation that $v(t_j)$ is locally bounded in $C^{2, \alpha}$ for all $j$ near any limiting point in $X\setminus D$, where $D$ is the singular set of $X$. After passing to a subsequence, we can assume that $(S^2, g(t_j))$ converges in Gromov-Hausdorff topology to $(S^2, d)$ equipped with a smooth Riemannian metric $g_\infty$ on $S^2\setminus D$, where $D$ is the singular set of $(S^2, d)$. Then $v(t_j)$ converges smoothly on $S^2\setminus D$ to $v_\infty$ satisfying on $X\setminus D$
\begin{equation}\label{soleq}
R(g_\infty) =
{1\over 2}\chi(S^2,\beta) + \Delta_\infty v_\infty.
\end{equation}
Furthermore, from the curvature bounds and injectivity radius bounds, for any domain $\mathcal{K}\subset\subset S^2\setminus D$, we can apply the local version of Hamilton's compactness theorem for the Ricci flow, i.e., there exist  domains $\mathcal{K}_j\subset\subset S^2\setminus D$ and diffeomorphisms $\Phi_j: \mathcal{K} \rightarrow \mathcal{K}_j$, such that the Ricci flow $g(t_j+t)$ for $t\in [0, 1]$ converges to a smooth family of Riemannian metrics $g_\infty(t)$ for $t\in [0, 1]$ on $\mathcal{K}$ satisfying the Ricci flow
$$\ddt{g_\infty(t)} = - Ric(g_\infty(t)) + \frac{1}{2} \chi(S^2, \beta) g_\infty(t), ~ g_\infty(0)=g_\infty. $$
In particular, $v(t_j+t)$ converges to $v_\infty(t)$ smoothly on $\mathcal{K} \times [0, 1]$ after reparametrization.

We claim that on $S^2\setminus D$, we have $$\nabla_\infty^2  v_\infty = \frac{1}{2} (\Delta v_\infty ) g_\infty.$$
Otherwise, there exists a domain $\mathcal{K}\subset\subset S^2\setminus D$ such that
$$\inf_{\mathcal{K}}|\nabla_\infty^2  v_\infty - \frac{1}{2} (\Delta v_\infty ) g_\infty|^2_{g_\infty}>0.$$
Then there exists some $\delta\in (0, 1)$ such that $$\inf_{\mathcal{K}\times [0, \delta]}|\nabla_{g_\infty(t)}^2  v_\infty (t)- \frac{1}{2} (\Delta_{g_\infty(t)} v_\infty(t) ) g_\infty(t)|^2_{g_\infty(t)}>0,$$
in particular,
$$\int_0^1\int_{\mathcal{K}} |\nabla_{g_\infty(t)} \nabla_{g_\infty(t)} v_\infty (t) -\frac{1}{2} (\Delta_\infty v_\infty(t) )g_\infty(t)|_{g_\infty(t)} dg_\infty(t)>0.$$
From the smooth convergence of the Ricci flow $g(t_j+t)$, this then implies that
$$\liminf_{j \rightarrow \infty}\int_{t=t_j}^{t_j+1} \int_{S^2} |\nabla^2 v - \frac{1}{2} (\Delta v) g|^2_g ~dg  ~dt>0.$$
Contradiction by Lemma \ref{entropy1}.

Hence $(g_\infty, v_\infty)$ is a shrinking gradient soliton on $X\setminus D$. In particular, $X_\infty =  \uparrow \dbar v_\infty$ is a holomorphic vector field on $X\setminus D$. From the partial $C^0$-estimate, $g_\infty$ extends to a K\"ahler current with bounded local potentials. Since $v_\infty$ is bounded in $W^{1, 2}(S^2)$ with respect to a fixed smooth metric on $S^2$, $X_\infty$ must extend to a holomorphic vector field on $S^2$.

\smallskip

We consider the following two cases:

\begin{enumerate}

\item $X_\infty$ is trivial.

\smallskip

In this case, $v_\infty$ is a constant and  the limiting metric is a constant curvature metric on $S^2\setminus D$.
Suppose that $D=\{P_1, ..., P_l\}$. We choose  holomorphic sections $\sigma_i \in H^0(S^2, -K_{S^2})$ such that  $\sigma_i$ vanishes at $P_i$ of order $2$. Let $\omega_{FS}\in c_1(S^2)$ be the standard smooth sphere metric on $S^2$. From the partial $C^0$ estimate and the fact that the limiting metric has constant curvature $1- \frac{1}{2}\sum_{i=1}^k\beta_i$  on $X\setminus D$, the limiting equation must be of the following form
$$\omega_\infty = \omega_{FS}  + \ddbar \varphi_\infty= e^{-(1- \frac{1}{2}\sum_{i=1}^k\beta_i)\varphi_\infty} (\omega_{FS})^{ 1- \frac{1}{2}\sum_{i=1}^k\beta_i }
 |\sigma_1|^{- \sum_{i=1}^k\beta_i} e^F$$
for some bounded potential $\varphi_\infty$ and some smooth pluriharmonic function $F$ on $S^2\setminus D$. Lemma \ref{onevar} implies that  the preceding equation can be rewritten in the following form
$$
\omega_\infty = \omega_{FS}  + \ddbar \varphi_\infty= \frac{e^{-(1- \frac{1}{2}\sum_{i=1}^k\beta_i )\varphi_\infty}}{\prod_{i=1}^l |\sigma_i|_{\omega_{FS}} ^{2\gamma_i}}  \omega_{FS} .
$$
This implies that $\omega_\infty$ is a conical constant curvature metric. The fact that the cone angle of each conical point is less than $2\pi$
implies that $\gamma_i\in (0,1/2)$, and
$\sum_{i=1}^l  2\gamma_i = \sum_{i=1}^k \beta_i$, from the Gauss-Bonnet formula.

\item $X_\infty$ is nontrivial.

\smallskip

Each nontrivial holomorphic vector field on $S^2$ can vanish at at most two distinct points and the imaginary part of $X$ is  a Killing vector field induced from an $S^1$-action. This implies that $D$ can have at most two points fixed by $X_\infty$ and the limiting soliton metric $g_\infty$ must be  rotationally symmetric. By the same argument as for the case of $X_\infty=0$ or directly by solving an ODE equation, we see that the limiting metric must be a conical shrinking gradient Ricci soliton metric on $S^2$ with $0$, $1$ or $2$ conical points. We denote $(S^2, \beta_\infty)$ the limiting conical sphere. In particular,  $\beta_\infty \neq 0$, by the Gauss-Bonnet formula.

\end{enumerate}

Combining the above, Lemma \ref{L1} is proved.

\bigskip

\subsection{The semi-stable case: $2\beta_{max} = \sum_{i=1}^k \beta_i$}

The goal of this section is to obtain a sequence converging to a conical constant curvature metric space along the conical Ricci flow on a semi-stable pair $(S^2, \beta)$.

The first step in the proof is to establish the following lower bound for the conical functional $F_\beta(\varphi)$: let $g(t)= g_\beta + \ddbar \varphi(t)$ be the solution of Ricci flow (\ref{ric}).  Then for any $\epsilon>0$, there exists $C_\epsilon>0$ such that for all $t\in [0, \infty)$ and $\varphi=\varphi(t)$,
\begin{equation}
\label{prop2}
F_{g_\beta}(\varphi) \geq - \epsilon \frac{\sqrt{-1}}{2\pi} \int_{S^2} \partial \varphi \wedge \dbar \varphi - C_\epsilon.
\end{equation}

\smallskip
To do this, we introduce for each $\epsilon>0$ the following approximation $F_{\beta,\epsilon}(\varphi)$ of the functional $F_\beta(\varphi)$ of (\ref{Fbeta}),
$$
F_{g_\beta, \epsilon}(\varphi) =  \frac{ \sqrt{-1}}{8\pi}\int_{S^2} \partial \varphi \wedge \dbar \varphi -   \frac{1}{2} \int_{S^2} \varphi dg_\beta-
{1\over
{1\over 2}\chi(S^2,\beta)- \epsilon} \log \int_{S^2}
e^{ -({1\over 2}\chi(S^2,\beta)- \epsilon) \varphi +h_\beta } dg_\beta
$$
We claim that for any $\epsilon\in (0,{1\over 2}\chi(S^2,\beta))$, the functional $F_{\beta,\epsilon}(\varphi)$ is bounded from below,
\begin{equation}
F_{\beta, \epsilon}(\varphi)
\geq -C_\epsilon,
\end{equation}
for all $\varphi \in PSH(S^2, g_\beta) \cap C^{2, \alpha}(S^2, \beta).$ This is because the Euler-Lagrange equation for the functional $F_{g_\beta,\epsilon}(\varphi)$ is a Monge-Amp\`ere equation which can be solved by the method of continuity for these values of $\epsilon$ (see \cite{T1} and \cite{BBGZ}). By the results of \cite{BBGZ}, the corresponding functional $F_{g_\beta,\epsilon}(\varphi)$ must be bounded from below by a constant for these values of $\epsilon$.

The lower bound for the functional $F_{\beta,\epsilon}$ implies the following lower bound for the functional $F_\beta$,
\begin{eqnarray}
\label{FFepsilon}
F_\beta(\varphi)
&\geq&
{\sqrt{-1}\over 8\pi}\int_{S^2}\partial\varphi\wedge\bar\partial\varphi
-
{1\over 2}\int_{S^2}\varphi dg_\beta
-
{2\over\chi(S^2,\beta)}
(\log\int_{S^2}e^{-({1\over 2}\chi(S^2,\beta)-\epsilon)\varphi}dg_\beta
-\epsilon {\rm inf}_{S^2}\varphi)
\nonumber\\
&\geq&
{\chi(S^2,\beta)-2\epsilon
\over
\chi(S^2,\beta)}F_{\beta,\epsilon}(\varphi)
+
{2\epsilon\over\chi(S^2,\beta)}
({\sqrt{-1}\over 8\pi}\int_{S^2}\partial\varphi\wedge\bar\partial\varphi
-
{1\over 2}\int_{S^2}\varphi dg_\beta
+
{\rm inf}_{S^2}\varphi)
\nonumber\\
&\geq&
{2\epsilon\over\chi(S^2,\beta)}
({\sqrt{-1}\over 8\pi}\int_{S^2}\partial\varphi\wedge\bar\partial\varphi
-
{1\over 2}\int_{S^2}\varphi dg_\beta+{\rm inf}_{S^2}\varphi)
-C_\epsilon.
\end{eqnarray}

By Lemma \ref{pere}, we know that the curvature and the diameter along the conical Ricci flow are uniformly bounded. We let $\omega_\beta = \hat \omega + \ddbar \psi$ and $\omega(t) = \hat \omega+ \ddbar \psi + \ddbar \varphi(t)$, where $\hat\omega\in [\omega_\beta]$ is a smooth Kahler metric and $\psi$ is a fixed continuous function in $PSH(S^2, \hat\omega)$. Let $\omega_j$ be the approximating smooth Kahler metrics for $\omega_t$ for a fixed $t$ as in Lemma \ref{metrapp}. Then there exists $C>0$ such that $\omega_j = \hat\omega + \ddbar (\psi + \varphi_j)$ satisfy
$$Ric(\omega_j) \geq -C \omega_j, ~ diam_{\omega_j}(S^2) \leq C, ~|\varphi_j - \varphi(t)|_{L^\infty(X)} \rightarrow 0$$
for all $j>0$.
Then the Green's functions $G_j$ for $\omega_j$ are uniformly bounded below for all $j$ and so
\begin{eqnarray*}
-\inf_{S^2} (\psi+\varphi_j) &\leq & - \frac{1}{2}\int_{S^2} (\psi+\varphi_j ) (\hat\omega+ \ddbar (\psi+\varphi_j) + K \\
&\leq& \int_{S^2} \frac{\sqrt{-1}}{4\pi} \partial \varphi_j \wedge \dbar \varphi_j + \frac{1}{2}\int_{S^2} \varphi_j \omega_\beta + K''
\end{eqnarray*}
for fixed $K$, $K'>0 $ because $\psi \in PSH(S^2, \hat\omega)\cap L^\infty(S^2)$. Hence by letting $j\rightarrow \infty$, we have
\begin{equation}\label{jfun}
-\inf_{S^2} \varphi   \leq \int_{S^2} \frac{\sqrt{-1}}{4\pi} \partial \varphi \wedge \dbar \varphi + \frac{1}{2} \int_{S^2} \varphi \omega_\beta + K''
\end{equation}
since $\varphi_j$ converges to $\varphi$ in $L^\infty$, where $K''$ only depends on $K$ and not on $t$.
Substituting this inequality into (\ref{FFepsilon}) gives the desired inequality (after a renaming of $\epsilon$).

\bigskip
We remark that in the last step, we can avoid using the lower bound of the Green's function for the evolving metrics. It suffices to approximate the evolving metrics $g(t)$ and $\varphi$ by smooth metrics and potentials so that the estimate (\ref{jfun}) holds uniformly for the approximation.

\bigskip
The estimate (\ref{prop2}) is slightly weak, since the ideal bound for the $F$ functional should be a uniform bound from below by a constant.
However, (\ref{prop2}) suffices for our purpose, which is to show that the curvature converges to a constant:

\begin{lemma} \label{curcon} There exists a sequence $t_j\rightarrow \infty$ such that the scalar curvature $R(t_j)$ converges uniformly to $1- \beta_{max} $ along the conical Ricci flow.

\end{lemma}

\noindent
{\it Proof}.  First, a straightforward calculation shows that
$$
\ddt{} F_{g_\beta} (\varphi_t) =  - \int_{S^2}v (1- e^{-v}) dg \leq 0,
$$
where $v $ is the Ricci potential defined in Lemma \ref{pere}.
Next, we claim that
\begin{equation}
\label{intv}
\inf_{t\in[T, \infty)}  \int_{S^2} v (1- e^{-v}) dg =0
\end{equation}
for any $T\geq 0$.
We prove the claim by contradiction. If not, then $\inf_{ [T' , \infty)}  \int_{S^2} v (1- e^{-v}) dg >\delta$ for some fixed $\delta>0$ and some sufficiently large $T'>0$. This implies that
$$F_{g_\beta}(\varphi_t) \leq - \delta t +  C_1,$$
for some fixed $C_1>0$.
On the other hand, by the estimate for $v$ established in Lemma \ref{pere}
and the fact that $\dot\varphi(t)=v$, $\varphi(t)$ has at worst linear growth in $t$ modulo a bounded time-dependent constant. Thus
$$\int_{S^2} \sqrt{-1}  \partial \varphi(t) \wedge \dbar \varphi(t) = \int_{S^2} \varphi(t) (\omega_\beta - \omega(t)) \leq At + C_2 $$
for some fixed $A, C_2>0$, where $\omega_\beta$ and $\omega(t)$ are the K\"ahler forms associated to $g_\beta$ and $\omega(t)$.
Therefore,
\begin{eqnarray*}
F_{g_\beta} (\varphi_t) + \epsilon \int_{S^2} \sqrt{-1}  \partial \varphi_t \wedge \dbar \varphi_t  &\leq& -\delta t + \epsilon \int_{S^2} \sqrt{-1} \partial \varphi_t \wedge \dbar \varphi_t+ C_1 \\
&\leq &- (\delta - A\epsilon) t + C_3 \rightarrow  - \infty
\end{eqnarray*}
 as $t\rightarrow \infty$ if we choose $\epsilon>0$ sufficiently small. This would contradict the estimate \ref{prop2}.

Since $|\nabla v |$ is uniformly bounded for $t\in [0, \infty)$,
and we have  $\kappa$-nonlocal  collapsing along the flow,
the inequality (\ref{intv}) implies that there exists a sequence $t_j\rightarrow \infty$ with
$\lim_{j\rightarrow \infty} \sup_{S^2} |v(t_j)| =0$. Furthermore, we have
$$\lim_{j\rightarrow \infty} \sup_{(z, t)\in S^2\times [t_j, t_j +1] } |v(z, t)| =0$$
since $\ddt{v} = \Delta v + v$ and $v$ is uniformly bounded.

Finally, one can apply the smoothing techniques in \cite{PSSW} in the conical setting and show that there exist a sequence $t_j' \in [t_j, t_j+1] \rightarrow \infty $ such that
$$ \lim_{j \rightarrow \infty} \sup_{S^2} ( |\nabla v(t_j')| + |R(t_j') -
{1\over 2}\chi(S^2,\beta) | )= 0.$$

This requires only the application of the maximum principle over $[t_j, t_j+1]$ combined with a family of barrier functions as in Lemma  \ref{maximum}, which is justified by the regularity of $v$. For example, for the bound on $|\nabla v|$, we would apply the maximum principle to $e^{-2t}(v^2+t|\nabla v|^2)-\epsilon \log |\sigma|^2_g$, and for $\Delta v$, to
$e^{-(t-1)}(|\nabla v|^2+ (t-1)\Delta v)-\epsilon |\sigma|_g^2)$, and let $\epsilon\to 0$. This completes the proof of Lemma \ref{curcon}.

\bigskip

By Lemma \ref{L1}, after passing to a subsequence, we can assume that $(S^2, g(t_j))$ in Lemma \ref{curcon} converges to a limiting space $(X, d)$.

\begin{lemma}\label{limpt} The limiting space $(X, d)$ has the following properties:

{\rm (1)}
The singular set $D$ consists of two points $p_\infty, q_\infty$ with weights $\beta_{max}$,

{\rm (2)} the conical Ricci flow converges in Gromov-Hausdorff topology to $( (S^2, \beta_\infty), g_\infty) $ with $\beta_\infty = \beta_{max} [p_\infty] + \beta_{max} [q_\infty ]$, and $g_\infty \in c_1(S^2)$ is the unique conical metric with constant curvature $1- \beta_{max}$,

{\rm (3)} the convergence is in $C^\infty$ on $S^2\setminus \{ p_\infty, q_\infty\}$.

\end{lemma}

\noindent
{\it Proof}.  By Lemma \ref{L1} and Lemma \ref{curcon}, $(S^2, g(t_j))$ converges to a conical metric of constant curvature $1-\beta_{max}$. The angle of $p_k$ is $2(1-\beta_{max})\pi$, which is the smallest angle. The point $p_k$ will converge to a limiting point $p_\infty$ in the limiting space $(X, d)$ along any convergent subsequence in the Gromov-Hausdorff topology. Applying the convergence results from the theory of Cheeger-Colding and volume comparison, we obtain
\begin{eqnarray*}
\frac{ Vol(B_d(p_\infty, r))}{Vol(B^{1-\beta_{max}}(r))} &=& \lim_{j\rightarrow \infty} \frac{ Vol(B_{g(t_j)} (p_k, r))}{Vol(B^{1-\beta_{max}}(r))} \\
&\leq & \lim_{j\rightarrow \infty} \lim_{r\rightarrow 0}    \left( \frac{ Vol(B_{g(t_j)} (p_k, r))}{Vol(B^{\inf_{S^2} R(g(t_j))}(r))}  \right) \left( \frac{Vol(B^{\inf_{S^2} R(g(t_j))}(r))}{Vol(B^{1-\beta_{max}}(r))}    \right) \\
&=& (1-\beta_{max}),
\end{eqnarray*}
because the curvature tends to $1-\beta_{max}$ uniformly as $j\rightarrow\infty$,  where $B^H(r)$ is the metric ball of radius $r$ on $S^2$ of constant curvature $H$ for $H>0$.
Therefore $p_\infty$ must be a singular point on $X$ by volume comparison. In particular, $p_\infty$ is a conical point with cone angle at most $2(1-\beta_{max}) \pi$.

\smallskip Applying Troyanov's stability condition for the existence of constant  curvature combined with $R= 1-\beta_{max}  $ on $S^2$, we can conclude that there can be only another conical point $q_\infty$ with the same cone angle as $p_\infty$, otherwise, the curvature must be strictly less than $1-  \beta_{max}$.

\smallskip Since the constant curvature metric in $c_1(S^2)$ with two conical points of cone angle $2\pi (1-\beta_k)$ is unique, the flow must converge to the same limiting space for any convergent subsequence.  This completes the proof of Lemma \ref{limpt}.

\bigskip

By the uniqueness of conical constant curvature metric on $(S^2, \beta_\infty)$, the limiting metric $g_\infty$ must be rotationally symmetric. We now want to show that the conical points $p_1, p_2, ..., p_{k-1}$ will merge into one point in the limiting space. This would complete the proof of Part 2 of Theorem \ref{main}.

\begin{lemma} \label{ptcon} Let $\mathcal{A}=\{ p_1, p_2, ..., p_{k-1}\}$. Then the diameter of $\mathcal{A}$ with respect to $g(t_j)$ converges to $0$ as $t_j\rightarrow \infty$. Furthermore, $\mathcal{A}$ converges to a conical point in the limiting space.

\end{lemma}

\noindent
{\it Proof}. First we pick $p_k$ and let $p_\infty$ be the limiting point of $p_k$ along the flow. It suffices to show that $\liminf_{j\rightarrow \infty} dist_{g(t_j)} ( p_k, \mathcal{A}) > 0$ by Lemma \ref{limpt}. This is because there is only one conical point $q_\infty $ other than $p_\infty$, and the limit of each $p_i$ must be a singular point by the volume comparison, hence $\mathcal{A}$ must converge to $q_\infty$.

\smallskip

 We will prove the proposition by contradiction. Suppose that a subset $\mathcal{A}'$ of  $\mathcal{A}$ converges to $p_\infty$ instead of $q_\infty$, say $q_1, ..., q_l$ with weights $\beta_{q_1}$, ..., $\beta_{q_l}$. We know that the limiting space is a football of constant curvature metric. Let $2L$ be the distance from $p_\infty$ and $q_\infty$ on $(X, d)$. Then $B_{g(t_j)}(p_k, L)$ converges to the half football $B_{g_\infty}(p_\infty, L)$ in Gromov-Hausdorff topology, furthermore, the convergence is smooth on $B_{g_\infty}(p_\infty, L)\setminus B_{g_\infty}(p_\infty, L/2)$.

\smallskip

 Let   $K = S^2 \setminus \{ \overline{B_{g_\infty} (p_\infty, L/4) } \cup \overline{ B_{g_\infty}(q_\infty, L/4)} \}$. Then for any $\epsilon>0$, there exists $T>0$ such that for $t_j>T$, there exists a diffeomorphism
$$\sigma_{t_j}: K  \rightarrow K(t_j) \subset S^2\setminus \beta$$
 such that
 $$ \| \sigma_{t_j}^*g(t_j) - g_\infty \|_{C^2(K, g_\infty)} < \epsilon. $$
Let $\eta$ be a smooth cut-off function on $(S^2\setminus \beta_\infty, g_\infty)$ such that
$0\leq \eta(z) \leq 1 $ on $S^2$ with $\eta =0$ on $S^2 \setminus B_{g_\infty}(q_\infty, L/2)$ and $\eta=1$ on $S^2 \setminus B_{g_\infty}(p_\infty, L/2)$.
Let $\tilde\sigma_{t_j}$ be a smooth diffeomorphism of $S^2$ which is a smooth extension of $\sigma_{t_j}$ to $S^2 \setminus \overline{B_{g_\infty}(q_\infty, L/4)}$. We then define a conical metric $\tilde g(t)$ by
\begin{equation}
\tilde g(t_j) = \eta g_\infty + (1-\eta) ~\tilde \sigma_{t_j}^* g(t_j).
\end{equation}
Obviously, $\tilde g(t_j) = g_\infty$ on $B_{g_\infty}(q_\infty, d/2)$ and $\tilde g(t_j) = \tilde\sigma_{t_j}^* g(t_j)$ on $B_{g_\infty}(p_\infty, d/2)$.

\smallskip Since $g(t_j)$ converges to $g_\infty$ on $K$, $\tilde g(t_j)$ converges to $g_\infty$ on $K$ smoothly as $t_j \rightarrow \infty$. This implies that $R(\tilde g(t_j))$ converge to $1- {1\over 2} \sum_{i=1}^k \beta_i$ in $L^\infty(S^2)$, and the total volume of $\tilde g(t_j)$ converges to 2, i.e.,
$\lim_{t\rightarrow \infty} \int_{S^2} d \tilde g(t_j) = 2$ because $B_{g(t_j)}(p_k, d/2)$ converges to $B_{g_\infty}(p_\infty, d/2)$ in Gromov-Hausdorff topology as well as in measure.   Therefore
\begin{equation}\label{limcov}
\lim_{t\rightarrow \infty} \int_{S^2} R(\tilde g(t_j)) d\tilde g(t) = 2- \sum_{i=1}^k \beta_i.
\end{equation}
On the other hand, by Gauss-Bonnet formula,
\begin{equation}\label{limcov2}  \int_{S^2} R(\tilde g(t_j)) d \tilde g(t) = 2- \beta_k - \sum_{i=1}^l \beta_{q_i} > 2- \sum_{i=1}^k \beta_i.
\end{equation}
Equations (\ref{limcov}) and (\ref{limcov2}) lead to contradiction by choosing $t$ sufficiently large.

\bigskip

This lemma illustrates why one cannot apply a local version of Hamilton's compactness theorem to the local $C^\infty$-convergence as in Proposition 5.3 in \cite{MRS}. This is because from the partial $C^0$-estimates, the gauge transformations come from the $\mathbb{C}^*$-action, and all  points but $p_k$ will converge to a single limiting conical point. Thus the injectivity radius will always tend to $0$ for generic points on $S^2\setminus \beta$.

\subsection{The unstable case: $2\beta_{max} > \sum_{i=1}^k \beta_i$}

In this section, we show that if $(S^2, \beta)$ is unstable, then any sequential limit cannot be a conical constant curvature metric space.

\begin{lemma} Suppose $(S^2, \beta_\infty )$ is a sequential limit of the conical Ricci flow (\ref{ricci}) on an unstable pair $(S^2, \beta)$. Then $\beta_\infty= \beta_{p_\infty} [p_\infty] + \beta_{q_\infty} [q_\infty]$ with $ 0\leq \beta_{q_\infty} < \beta_{q_\infty} <1$. Therefore, the limiting soliton metric can not be a constant curvature metric.

\end{lemma}

\noindent
{\it Proof}. Suppose $(S^2, g(t_j))$ converges to a limiting conical shrinking gradient Ricci soliton $((S^2, \beta_\infty), g_\infty)$. By the volume comparison, $p_k$ must converge to a limiting conical point, say $p_\infty$, such that the cone angle of $g_\infty$ at $p_\infty$ must be at most $2\pi(1-\beta_k)$ by volume comparison as the curvature is uniformly bounded.  On the other hand, by the boundedness of $R$ and smooth convergence of $R$ on $S^2\setminus D$, we have
 $$\int_{S^2} R(g_\infty) d g_\infty = 2- \sum_{i=1}^k \beta_i. $$
Suppose $\beta_\infty = \beta_{p_\infty} [p_\infty] + \sum_{i=1}^l \beta_{q_i} [q_i] $ with distinct points $p_\infty$, $q_1, ..., q_l$.
 Then $\beta_{p_\infty} \geq \beta_k$  and by  the Gauss-Bonnet formula, we have
$$ \beta_{p_\infty} + \sum_{i=1}^l \beta_{q_i}= \sum_{i=1}^k \beta_k$$
and so
$$ \sum_{i=1}^l \beta_{q_i} \leq \sum_{i=1}^{k-1} \beta_i < \beta_k \leq \beta_{p_\infty}. $$
This contradicts Troyanov's stability condition and so $(S^2, \beta_\infty)$ does not admit a conical constant curvature metric.


\subsection{Uniqueness and uniform convergence}

In this section, we shall complete the proof of Theorem \ref{main}. The major step is to pass from the sequential convergence established in the previous section to a full convergence. We will always assume that $(S^2, \beta)$ is either semi-stable or unstable.

\smallskip

 \begin{lemma} \label{5solim}  Assume that $(S^2, \beta_\infty, g_{sol, \beta_\infty})$ is the limit of a sequence $(S^2, \beta, g(t_\mu))$ along the conical Ricci flow as $\mu \rightarrow \infty$.
Then $(S^2, \beta_\infty, g_{sol, \beta_\infty})$ is a rotationally symmetric conical gradient shrinking Ricci soliton. Furthermore,  $$\beta_{p_\infty} = \sum_{i\in I} \beta_i, ~~\beta_{q_\infty} = \sum_{j\in J} \beta_j, $$
for some $ I \subset \{ 1, 2, ..., k\}$ and $J = \{1, 2, ..., k\} \setminus I. $

 \end{lemma}

\noindent
{\it Proof}. If $g_{sol, \beta_\infty}$ has constant curvature, then by Lemma \ref{limpt} using Troyanov's stability condition, $(S^2, \beta_\infty)$ must be semi-stable and $$\beta_\infty = \beta_{max} [p_\infty] + \beta_{max} [q_\infty] $$ for two distinct points $p_\infty$ and $q_\infty$. Therefore, the limiting metric must be always rotationally symmetric.

Without loss of generality, we can assume that for some $I \subset \{ 1, 2, ..., k\}$ and $J = \{1, 2, ..., k\} \setminus I$, $$p_i \rightarrow p_\infty, ~ p_j \rightarrow q_\infty$$
 for $i\in I$ and $j\in J$. It suffices to show that the weights at $p_\infty $ and $q_\infty$ satisfy
 $$\beta_{p_\infty} = \sum_{i\in I} \beta_i, ~~\beta_{q_\infty} = \sum_{j\in J} \beta_j. $$
 This can be shown by the same arguments in the proof of Lemma \ref{ptcon} by gluing and by the Gauss-Bonnet formula, because the curvature is uniformly bounded and converges uniformly away from $p_\infty$ and $q_\infty$.

\begin{lemma} \label{finitetp} Let $\mathcal{S}$ be the set of all  conical shrinking gradient Ricci solitons $(S^2, \beta_\infty, g_{sol, \beta_\infty})$ which arise as sequential limits   for the Ricci flow. Then $\mathcal{S}$ is a finite set.
\end{lemma}

\noindent
{\it Proof}. For fixed $\beta_\infty$, the conical shrinking soliton $(S^2, \beta_\infty, g_{sol, \beta_\infty})$ is unique. The corollary immediately follows from Lemma \ref{5solim} and the fact that there are only finitely many combinations of $I \sqcup J = \{1, 2, ..., k\}$.

\begin{lemma} \label{uniqe} Let $g(t)$ be the solution of the Ricci flow. Then $(S^2, \beta, g(t))$ converges uniformly in Gromov-Hausdorff topology to a shrinking gradient conical Ricci soliton $(S^2, \beta_\infty, g_{sol, \beta_\infty})$ for $t\rightarrow \infty$.

\end{lemma}

\noindent
{\it Proof}.  The proof is by contradiction. Suppose not. Then  there exist two sequences of solutions for the conical Ricci flow $g(t_l)$ and $g(t_l')$ converging to two distinct conical shrinking solitons $(S^2, \beta', g_{sol}')$ and $(S^2, \beta'', g_{sol}'')$ in Gromov-Hausdorff topology, as $l\rightarrow \infty$. In particular there exists $D>0$ such that
$$d_{GH}((S^2, g'), (S^2, g'')) =2D$$
and $L>0$ such that for all $l>L$,
$$d_{GH}((S^2, g(t_l)), (S^2, g(t_l'))) >D.$$
Without loss of generality, we can assume $t_l'> t_l$ for each $l$. Then we consider the function
$$f_l(s) = d_{GH} ((S^2, g((1-s)t_l + s~ t_l')) , (S^2, g(t_l))), ~ s\in [0, 1]. $$
Since the conformal factor of $g(t)$ with respect to $g(0)$ is continuous, $f(s)$ is a continuous function with $f_l(0)=0$ and $f_l(1)>D$.

Therefore, for any $d\in [0, D]$, there exist a sequence $g(t_{l, d})$ such that $$d_{GH}((S^2, g( t_{l,d} ) ) , (S^2, g(t_l))) =d. $$
After passing to a subsequence, $(S^2, g(t_{l, d}))$ converges to a conical shrinking soliton $(S^2, \beta_d, g_{sol, d})$ satisfying
$$d_{GH} ((S^2, g_{sol, d} ), (S^2, g_{sol}') )=d. $$
This implies that there are infinitely many distinct limiting conical shrinking solitons from the conical Ricci flow. This contradicts Lemma \ref{finitetp}.

\bigskip

We would also like to exclude the possibility that even though the limiting soliton is unique, the limiting conical points may be limits of different sets of conical points along the flow.

\begin{lemma} Let $g(t)$ be the solution of the Ricci flow on $(S^2,\beta)$ with a regular initial metric. Then tre exists a unique combination $I\sqcup J =\{1, 2, ..., k\}$, possibly dependent on the initial metric,  such that
$$p_i \rightarrow p_\infty, ~ p_j \rightarrow q_\infty$$
as $t\rightarrow \infty$, for any $i\in I$ and $j\in J$.

\end{lemma}

\noindent
{\it Proof}. Assume otherwise. Then there exists $p_i$ such that along a sequence $g(t_l)$, $p_i \rightarrow p_\infty$ and along another sequence $g(t_l')$, $p_i \rightarrow q_\infty$. Let $(S^2, g_{sol})$ be the limiting shrinking Ricci soliton. Then there must exist $\epsilon>0$ and a sequence $t_l''$ such that along $g(t_l'')$, %
$$d_{(S^2, g(t_l''))\sqcup (S^2, g_{sol})}( p_i, p_\infty) > \epsilon, ~ d_{(S^2, g(t_l''))\sqcup (S^2, g_{sol})}( p_i, q_\infty) > \epsilon $$
by continuity. This is a contradiction, because $p_i$ must converge to either $p_\infty$ and $q_\infty$ along $(S^2, g(t_l''))$ as $l\rightarrow \infty$.

\medskip

Combining the above results, we can now finish the proof of Theorem \ref{main}.

\medskip
It remains to prove the limiting soliton does not depend on the choice of the initial metric. Assume that it does. We discuss the following two possible cases.

\smallskip
In the first case, there exist two initial conical metrics $g_0$ and $g_0'$ such that the solutions $g(t)$ and $g'(t)$ along the conical Ricci flow converge to two distinct limiting solitons $(S^2, \beta_\infty, g_{sol})$ and $(S^2, \beta_\infty', g_{sol}')$. Suppose
 $$g_0' = e^{u_0} g_0. $$
 We consider a family of conical metrics
 $$g_{0, s} = e^{su_0} g_0, ~ s\in [0, 1]. $$
 Then $g_{0,s}$ satisfies the assumptions in Theorem \ref{main} for all $s\in [0, 1]$. Then we have a family of solutions $g_s(t) =e^{u_s(t)} g_0$ of the conical Ricci flow (\ref{ricu}) for $s\in [0, 1]$ with $g_s(0)=g_{0,s}$. By Yin's estimates, $u_s(t)$ is both continuous in $t$ and $s$ because $\dot u_s = \frac{\partial u_s}{\partial s}$ satisfies
 $$\ddt{\dot u_s} = e^{-u_s} \Delta_{g_\beta} \dot u_s  - e^{-u_s} \dot u_s \Delta_{g_\beta} u_s + e^{-u_s} \dot u_s R_0, ~ \dot u_s|_{t=0} = u_0. $$
 This implies that $(S^2, \beta, g_s(t))$ is continuous in both $t$ and $s$ in Gromov-Hausdorff topology. Then by the same arguments in the proof of Lemma \ref{uniqe}, there exist infinitely many distinct limiting shrinking solitons by varying $s$ in $g_s(t)$. This contradicts Lemma \ref{finitetp}.

 \smallskip
 In the other case, there exist $g_0$ and $g_0'$ such that the solutions $g(t)$ and $g'(t)$ along the conical Ricci flow with $g(0)=g_0$ and $g'(0)=g'_0$ converge to the same $(S^2, \beta_\infty, g_{sol})$ but there exists $d>0$ such that
 $$ \limsup_{t\rightarrow \infty} d_{GH}( (S^2, g(t)), (S^2, g'(t))) > d. $$
 By the same previous argument of constructing a family of Ricci flow solutions with parameter $s\in [0, 1]$, one can show there exist infinitely many distinct limiting shrinking Ricci soliton. That is again a contradiction.


 \section{Proof of Theorem \ref{main2} }

 In this section, we shall establish Theorem \ref{main2}. We conjecture that if $(S^2, \beta)$ is unstable, then the limiting soliton must be $(S^2, \beta_\infty)$ with
$$\beta_\infty = \beta_{p_\infty} [p_\infty] + \beta_{q_\infty} [q_\infty] , ~ \beta_{p_\infty} = \beta_k. $$
 We show now that, under some additional condition, we do obtain
 the uniform convergence of the Ricci flow towards such a soliton.  We define the following normalized $W$-functional
 \begin{equation}
W(g, f) =  \int_{S^2}   \left( \frac{1}{2- \sum_{i=1}^k \beta_i} (R + |\nabla f|^2) + f  \right) e^{-f} dg, ~ \int_{S^2} e^{-f} dg =2, ~ g\in c_1(S^2).
 \end{equation}
 Note that the singular time for the unnormalized conical Ricci flow on $(S^2, \beta)$ is
 $$\tau = (2-\sum_{i=1}^k\beta_i)^{-1}$$ and so
 $$W(g, f) = 2 W\left(g, f+ \log \frac{2-\sum_{i=1}^k\beta_i}{2\pi}, (2-\sum_{i=1}^k\beta_i)^{-1}) \right) +4 - 2 \log \frac{2-\sum_{i=1}^k\beta_i}{2\pi}   $$
 with
 $$  (4\pi\tau)^{-1} \int_{S^2} e^{-\left(f+ \log \frac{2-\sum_{i=1}^k\beta_i}{2\pi}\right)}  = 1.$$
As usual, we define
 $$   \mu (g) = \inf_f  W (g, f), ~ \int_{S^2} e^{-f} dg =2, $$
 where the infimum is taken over functions $f$ satisfying the combined conditions  (\ref{Fspace}).

 Suppose $((S^2, \beta), g_{sol})$ is a gradient shrinking soliton with $g_{sol}\in c_1(S^2)$. Then $g_{sol}$ is rotationally symmetric and satisfies
 \begin{equation}
 R(g_{sol}) = (1- \frac{1}{2}( \beta_p + \beta_q)) + \Delta_{g_{sol}} \theta_{sol},  ~  \nabla_{g_{sol}}^2 \theta_{sol} = \frac{1}{2} (\Delta_{g_{sol}} \theta_{sol} ) g_{sol}, ~\int_{S^2}  e^{\theta_{sol}}dg_{sol} =2
 \end{equation}
  for a unique $\theta_{sol}$.
  \begin{lemma} Suppose $((S^2, \beta), g_{sol})$ is a gradient shrinking soliton. Then
\begin{equation}
  W(g_{sol}, -\theta_{sol}) = 1  - \int_{S^2} \theta_{sol} e^{\theta_{sol}} dg_{sol}
\end{equation}

\end{lemma}

 \noindent
 {\it Proof}. Following the same argument by Hamilton \cite{H}, one can show that
 $$R- \Delta \theta_{sol} = (1-\frac{1}{2} \sum_{i=1}^k
 \beta_i), ~~R+ |\nabla \theta_{sol} |^2 = - (1- \frac{1}{2}\sum_{i=1}^k \beta_i) (\theta_{sol} +C). $$
Integrating by parts, we obtain
 $$C = 1+ \frac{1}{2}\int_{S^2} \theta_{sol} e^{\theta_{sol}} dg_{sol}$$
 and the lemma immediately follows.
 \bigskip

 We now compare $ W(g_{sol}, -\theta_{sol})$ for different markings,
 and establish a monotonicity formula for different conical shrinking soliton metrics.

 \begin{lemma} Let  $(S^2, \beta)$ and $(S^2, \beta')$ be two conical spheres with $\beta= \beta_p [p] + \beta_q[q]$, $\beta'= \beta'_{p'} [p'] + \beta'_{q'}[q']$, $\beta_p, \beta_q, \beta_{p'}, \beta_{q'}\in [0,1)$.  Let  $g_{sol, \beta}, g_{sol, \beta'} \in c_1(S^2)$ be the shrinking gradient soliton metrics on $(S^2, \beta)$ and $(S^2, \beta')$.
 If $\beta_p + \beta_q = \beta'_{p'}+\beta'_{q'}$ and $|\beta_p-\beta_q| < |\beta'_{p'}-\beta'_{q'}|$, then
 $$W(g_{sol, \beta}, -\theta_{sol, \beta})  > W(g_{sol, \beta'}, -\theta_{sol, \beta'} )$$

 \end{lemma}

 \noindent
 {\it Proof}.   It suffices to calculate the integral $\int_{S^2} \theta_{sol} e^{\theta_{sol}} dg_{sol}.$ We can apply the calculations in \cite{DGSW}, since the soliton metrics are toric.
 The polytope associated to $(S^2, c_1(S^2))$ is $P=[-1, 1]$ with defining functions $l_0(x) = 1-x\geq 0 $ and $l_\infty(x) = 1+x \geq 0$.  The soliton equation for $g_{sol}$ is given by
 $$Ric(g_{sol}) = (1-\frac{1}{2}( \beta_0 +\beta_\infty)) g_{sol} + L_\xi g_{sol} + \beta_0[D_0] + \beta_\infty[D_\infty],$$
 where $D_0$, $D_\infty$ are the two points fixed by the torus action, $\xi$ is a holomorphic vector field.  We let $\eta = |\beta_0-\beta_\infty|$.
 By Theorem 1.1 in \cite{DGSW}, one can solve the above equation if and only if
 $$\beta_0 [D_0] + \beta_\infty [D_\infty] = (1- (1-\frac{1}{2}(\beta_0+\beta_\infty) l_0(\tau) ) [D_0] + (1- (1-\frac{1}{2}(\beta_0+\beta_\infty) l_\infty (\tau) ) [D_\infty] .$$
 Immediately one has $$\tau = \frac{\beta_\infty - \beta_0} {2 - \beta_0-\beta_\infty}. $$
 Obviously, $|\tau|$ is an increasing function in $\eta$ since $\beta_0+\beta_\infty$ is fixed. Then one can uniquely solve $c$ from the following equation
 $$ \tau = \frac{\int_{-1}^1x e^{cx} dx}{\int_{-1}^1 e^{cx}dx}. $$
 In particular,
 $$\tau'(c) = \left(\int_{-1}^1 e^{cx} dx \right)^{-2}\left( \int_{-1}^1 x^2 e^{cx} dx \int_{-1}^1 e^{cx}dx - \left(\int_{-1}^1 x e^{cx}dx \right)^2 \right)>0.$$
 Therefore $|c|$ is an increasing function in $|\tau|$.
 From \cite{WZ, DGSW}, $$\theta_{sol} = \log \frac{2e^{cx}}{\int_{-1}^1 e^{cx}dx} $$ and
 $$F(c) =  \int_{S^2}  \theta_{sol} e^{\theta_{sol} }dg_{sol} =A^{-1} \int_{-1}^1(cx - \log A) e^{cx} dx, ~ A = \frac{1}{2} \int_{-1}^1 e^{cx} dx. $$
 Straightforward calculations show that
 $$F'(c) = \frac{c}{2A^2} \left( \int_{-1}^1x^2 e^{cx} dx  \int_{-1}^1 e^{cx}dx - \left(\int_{-1}^1 xe^{cx}dx \right)^2    \right). $$
 Therefore $F'(c)>0$ if $c>0$ and $F'(c)<0$ if $c<0$, and it immediately  implies that $ \int_{S^2}  \theta_{sol} e^{\theta_{sol} }dg_{sol} $ is strictly increasing in terms of $\eta$.
 This completes the proof of the lemma.

 \bigskip

Let $(S^2, \beta)$ be the sphere with
marked points $\beta= \sum_{i=1}^k \beta_i [p_i]$,  and $\beta_k > \sum_{i<k} \beta_i$. We let $I\sqcup J =\{1, 2, ..., k\}$ be a division of $\{1, 2, ..., k\}$  and define a
sphere with new marked points $(S^2, \beta_{I, J} )$   by
 $$\beta_{I, J}=  \sum_{i\in l} \beta_{i}[p] + \sum_{j\in J} \beta_{j} [q].$$
Then we can order the finite set $\{ \mu(g_{sol, \beta_{I, J}}, -\theta_{g_{sol, \beta_{I, J}}}) \}_{I, J} $ by
$$\mu_1> \mu_2\geq \mu_3\geq ... \geq \mu_N$$
for some $N$ and $\mu_1=    W(g_{sol, \beta_{I, J}}, -\theta_{g_{sol, \beta_{I, J}}}))$ with $I=\{k\}$ and $J=\{1, 2, ..., k-1\}$.

 \begin{lemma} \label{soon} Assume that the initial metric $g_0\in c_1(S^2)$ satisfies $  \mu (g_0) > \mu_2$. Let $(S^2, \beta_\infty)$ be a sphere with
 $$\beta_\infty = \beta_k[p_\infty] + (\sum_{i<k} \beta_i) [q_\infty].$$
 Then

 {\rm (1)}  the conical Ricci flow (\ref{ricci}) converges to the unique shrinking gradient soliton $((S^2, \beta_\infty), g_\infty)$ in Gromov-Hausdorff topology.

 \smallskip

 {\rm (2)} the convergence is smooth on $S^2\setminus \beta_\infty$,

 \smallskip

 {\rm (3)} $p_k$ converges to $p_\infty$ and $p_i$ converges to $q_\infty$ for all $i=1, ..., k-1$.

 \end{lemma}

 \noindent
 {\it Proof}.  Let $(S^2, \beta', g_{sol}')$ be the unique limiting shrinking soliton of the Ricci flow on $(S^2, \beta)$. It suffices to show that $\beta' = \beta_\infty$. Suppose not.   Since $g(t)$ converges to $g_{sol}'$ smoothly on any compact set $K$ of $S^2\setminus \beta'$, $v$ also converges to $\theta_{sol}$ smoothly on $K$. Suppose $\beta'= \beta_{p_\infty'} p_\infty' + \beta_{q_\infty'} q_\infty'$ and $R(g_{sol}') = (1-\frac{1}{2} \sum_{i=1}^k \beta_i) + \Delta_{g_{sol}'} \theta_{sol}'$ with $\int_{S^2} e^{\theta_{sol}'} dg_{sol}'=2$. For any $\epsilon>0$, there exist $\delta>0$ and $T>0$ such that
$$Vol_{g_{sol}'}(B_{g_{sol}'}(p_\infty', \delta) ) + Vol_{g_{sol}'}(B_{g_{sol}'}(q_\infty', \delta) ) \leq A^{-1} \epsilon, $$
where $A = 100 \sup_{S^2 \times [0, \infty)}  (  |R(g(t)|  +  |\nabla v|^2) e^{v} $,
and
on $K= S^2\setminus \left( B_{g_{sol}'}(p_\infty', \delta) \cup B_{g_{sol}'}(q_\infty', \delta) \right)$,
$$\| \sigma_t^* g(t) - g_{sol}'\|_{C^2(K, g_{sol}')} + \|v- \theta_{sol}'\|_{C^2(K, g_{sol}')} \leq 100^{-1} \epsilon, $$
for $t>T$, where $\sigma_t$ is a diffeomorphism in a neighborhood of $K$ at $t$. Then immediately we have,
$$ |   W(g(t), - v(t))  - W(g_{sol}', -\theta_{sol}') | < \epsilon   $$
for $t>T$.
Therefore,
	 $$\mu_2 <  \mu (g_0) \leq \lim_{ t\rightarrow \infty}   W (g(t), v(t)) =   W (g_{sol}', \theta_{sol}') \leq \mu_2.$$
This is a contradiction and so $\beta'=\beta_\infty$.  Lemma \ref{soon} is proved, and hence so is Theorem \ref{main2}.

\bigskip
\bigskip

\noindent {\bf{Acknowledgements:}} The second named author would like to thank Xiaochun Rong for many stimulating  discussions.

{\noindent \footnotesize $^*$ Department of Mathematics\\
Columbia University, New York, NY 10027\\

\noindent $\dagger$ Department of Mathematics\\
Rutgers University, Piscataway, NJ 08854\\

\noindent $\ddagger$ Department of Mathematics\\
Rutgers University, Newark, NJ 07102\\ }

\end{document}